\documentclass{article}

\usepackage{a4wide}
\usepackage{amsmath}
\usepackage{txfonts}
\usepackage{mathtools}
\usepackage{mdwlist}

\newtheorem{remark}{Remark}%\numberwithin{remark}{section}
\newtheorem{theorem}{Theorem}%\numberwithin{theorem}{section}

\newcommand{\pt}{\partial}
\newcommand{\RR}{\mathbb{R}}
\renewcommand{\phi}{\varphi}
\renewcommand{\rho}{\varrho}
\renewcommand{\theta}{\vartheta}
\newcommand{\scal}[2]{\left\langle#1,#2\right\rangle}
\newcommand{\m}[1]{\mathcal{#1}}
\newcommand{\E}{\mathrm{e}}
\newcommand{\D}{\mathrm{  d}}
\newcommand{\abs}[1]{\left|#1\right|}
\newcommand{\norm}[1]{\left\|#1\right\|}

\title{A review of maximum-norm a posteriori error bounds for
       time-semidiscretisations of parabolic equations}

\author{Torsten Lin\ss\thanks{Fakult\"at f\"ur Mathematik und Informatik,
        FernUniversit\"at in Hagen,
        Universit\"atsstra{\ss}e 11,
        58095 Hagen,
        Germany,
        \texttt{torsten.linss@fernuni-hagen.de},\ \
        \texttt{martin.ossadnik@fernuni-hagen.de}}
   \and Natalia Kopteva\thanks{Department of Mathematics and Statistics,
        University of Limerick, Limerick, V94 T9PX, Ireland,
        \texttt{natalia.kopteva@ul.ie}}
   \and Goran Radojev\thanks{Department of Mathematics and Computer Science, Faculty of Sciences,
        University of Novi Sad, Trg Dositeja Obradovi\'ca~4, 21000 Novi Sad,
        Serbia,
        \texttt{goran.radojev@dmi.uns.ac.rs}}
   \and Martin Ossadnik\setcounter{footnote}{0}\footnotemark{}
}

\begin{document}

\maketitle

\begin{abstract}
  A posteriori error estimates in the maximum norm are studied for
  various time-semidiscretisations applied to a class of linear parabolic
  equations.
  We summarise results from the literature and present some new improved
  error bounds. 
  Crucial ingredients are certain bounds in the $L_1$ norm for the Green’s
  function associated with the parabolic operator and its derivatives.



\end{abstract}

%%%%%%%%%%%%%%%%%%%%%%%%%%%%%%%%%%%%%%%%%%%%%%%%%%%%%%%%%%%%%%%%%%%%%%%%%%%%%%%%
\section{Introduction}

Consider the linear parabolic equation:
\begin{subequations}\label{problem}
\begin{alignat}{2}
 \m{K}u \coloneqq \pt_t u  + \m{L} u & = f\,,
            &\quad& \text{in} \quad Q\coloneqq\Omega \times (0,T],\\
 \intertext{with a second-order linear elliptic operator $\m{L}$
    in a spatial domain $\Omega\subset\RR^n$ with Lipschitz boundary and
    some function $f\colon [0,T]\to L_2(\Omega)$, subject to the
    initial condition}
   u(x,0) & = u^0(x)\,, && \text{for} \quad x\in\bar{\Omega}, \\
 \intertext{and the Dirichlet boundary condition}
   u(x,t) & = 0\,, && \text{for} \quad (x,t)\in \pt\Omega \times [0,T].
\end{alignat}
\end{subequations}

Following \cite{MR1335652} and \cite{MR2519598}, the authors of the present
study have published a number of results on residual-type a posteriori
error estimates in the maximum norm for parabolic equations utilising and
merging various approaches and considering various classes of temporal
discretisation~\cite{MR2519598, MR2629992, MR3032709, MR3056758, RadLin22,
MOssad22}.
In this survey, we review these results in a unified manner.
Reexamining those results and their proofs, we are able to present some
improvements, namely for the implicit Euler method, the Crank-Nicolson method
and the dG(1)-method.
These improvements are made possible by using local, time-slice wise bounds
for the Green's function, rather than global stability results.
Details will be highlighted in the course of the paper.
We also present some new results (most notably Theorems~\ref{theo:LR} and~\ref{theo:CN2}).
Furthermore, numerical results are given to compare the various approaches.

The general idea is to represent the error (at final time $T$) by means of the
Green's function and the residual.
To this end bounds in the $L_1$-norm on the Green’s function associated with $\m{K}$
and its time-derivatives are required, see~\S\ref{sect:green} for details.

In the present paper we study semidiscretisations in time only.
However, these are essential building blocks in deriving error estimates for
full (space-and-time) discretisations.
Using so called elliptic reconstructions, they can be combined with error
estimators for discretisations of elliptic problems to give error bounds for
parabolic problems.

The paper is organised as follows.
In Section~\ref{sect:green} we specify our general assumptions for the a
posteriori error analysis, in particular the stipulate the validity of 
certain bounds for the Green's function of the parabolic problem.
Thereafter, we present result for various discretisations:
\begin{itemize*}
  \item the simple first-order implicit Euler method (\S\ref{sect:euler}),
  \item the second-order Crank-Nicolson method (\S\ref{sect:CN}),
  \item an extrapolated Euler method of 2nd order (\S\ref{sect:extra}),
  \item the third-order discontinuous Galerkin-Radau method
        (dG(1), \S\ref{sect:dG1}), and finally
  \item the backward-differentiation formula of order $2$ (\S\ref{sect:BDF2}).
\end{itemize*}
We complement the theoretical finds with results of numerical experiments.
The test problem is introduced in \S\ref{sect:test}.

\textbf{Notation:} For functions $w\colon \Omega\times[0,T]\to\RR$ we shall use
the shortend notation $w(t)\coloneqq w(\cdot,t)$ which for each time $t$ is a
function mapping from $\Omega$ to $\RR$.

\section{The Green's function}
\label{sect:green}

In this section we consider the Green’s function G associated with the operator
$\m{K}$ in~\eqref{problem}.
It will be used to express the error of a numerical approximation in terms of
its residual in the differential equation.
For definitions and properties of fundamental solutions and Green’s functions of
parabolic operators, we refer the reader to the survey by Friedman
\cite{0144.34903}, in particular Chapter 1.

For fixed $x\in\Omega$, the Green's function associated with $\m{K}$
and $x$ solves
\begin{gather*}
  \pt_t\m{G}+\m{L}^*\m{G} = 0, \ \ \text{in} \ \Omega\times\RR^+, \ \
  \m{G}\bigr|_{\pt\Omega} = 0, \ \ \m{G}(0) = \delta_x=\delta(\cdot-x)\,,
\end{gather*}
with $\delta$ denoting the Dirac $\delta$-distribution.
Let $\scal{\cdot}{\cdot}$ denote both the duality pairing on
\mbox{$H^{-1}(\Omega)\times H_0^1(\Omega)$} and the $L_2(\Omega)$ scalar product.
Then for all \mbox{$w\in W^{1,1}\bigl([0,T],H_0^1(\Omega)\bigr)$} and
$t\in(0,T]$, we have
%%
%\begin{align*}
%  0 & = \int_0^t \scal{\pt_t\m{G}(t-s) + \m{L}^*\m{G}(t-s)}{w(s)}\D s \\
%    & = \scal{\m{G}(t)}{w(0)} - \scal{\m{G}(0)}{w(t)}
%           + \int_0^t \scal{\m{G}(t-s)}{\pt_t w(s)} \D s
%           + \int_0^t \scal{\m{G}(t-s)}{\bigl(\m{L} w\bigr)(s)} \D s \\
%    & = \scal{\m{G}(t)}{w(0)} - w(x,t)
%           + \int_0^t \scal{\m{G}(t-s)}{\bigl(\m{K} w\bigr)(s)} \D s\,.
%\end{align*}
%%
%This gives
%%
\begin{gather}\label{green-rep}
  w(x,t) = \scal{\m{G}(t)}{w(0)}
              + \int_0^t \scal{\m{G}(t-s)}{\bigl(\m{K}w\bigr)(s)} \D s.
\end{gather}
We will make frequent use of this representation of a function $w$ in terms
of its residual $\m{K}w$.

Throughout the paper we shall assume there exist non-negative constants
$\kappa_0$, $\kappa_1$, $\kappa_2$, $\kappa_1'$, $\kappa_2'$ and
$\gamma$ such that (with formally setting $\kappa_0'=0$)
\begin{gather}\label{source:ass}
  %\norm{\m{G}(t)}_{1,\Omega} \le \kappa_0\,\E^{-\gamma t} \eqqcolon \phi_0(t),
  %   \quad
  \norm{\pt_t^{p} \m{G}(t)}_{1,\Omega}
     \le \left(\frac{\kappa_{p}}{t^p} +\kappa_p'\right)
                \,\E^{-\gamma t} \eqqcolon \phi_p(t), \quad
     \text{for all} \ x\in\bar\Omega, \ t\in[0,T], \text{and} \  p=0,1,2.
\end{gather}
Here $\norm{\cdot}_{p,\Omega}$, \mbox{$p\in[1,\infty]$}, denotes the standard
norm in $L_p(\Omega)$.
A number of problems that satisfies these assumptions are gathered
in~\cite[\S2.1]{MR3720388}.
There results from various sources are summarised, including
\cite[\S2.2]{MR2519598}, \cite{MR1423289}
\cite[\S12]{MR3056758} and the case of a singularly perturbed problem
in~\cite[\S2]{MR3032709}.

The rest of this section is rather technical as we will precompute some coefficients
that feature in our error bounds later.
They appear after H{\"o}lder's inequality and~\eqref{source:ass} have been applied
to integrals involving (derivatives of) the Green's function.
Those integrals are of the form
\begin{gather*}
  \int_{t_{j-1}}^{t_j} \pi(s) \phi_p(T-s)\D s
     \ \ \ \text{with} \ \ \ 0\le t_{j-1}<t_j\le T,
     \ \ p=0,1,2, \ \ \text{and a function} \ \pi.
\end{gather*}
These are bounded as follows
\begin{gather}\label{intG_p}
  \abs{\int_{t_{j-1}}^{t_j} \pi(s) \phi_p (T-s)\D s}
     \le \E^{-\gamma(T-t_j)} 
     \int_{t_{j-1}}^{t_j} \abs{\pi(s)}
            \left(\frac{\kappa_{p}}{(T-s)^p} +\kappa_p'\right) \D s\,.
\end{gather}
For example,
\begin{gather}\label{intG}
  \abs{\int_{t_{j-1}}^{t_j} \pi(s) \phi_0(T-s) \D s}
     \le \kappa_0 \E^{-\gamma(T-t_j)}
       \int_{t_{j-1}}^{t_j} \abs{\pi(s)} \D s\,, \\
 \label{intG_t}
   0 \le \int_{t_{j-1}}^{t_j} \phi_1(T-s) \D s
    = \E^{-\gamma(T-t_j)} \theta_j, \quad
        \theta_j \coloneqq \left\{\kappa_1 \ln\left(1+\frac{\tau_j}{T-t_j}\right)
                             + \kappa_1'\tau_j\right\}
 \intertext{and}\label{intG_t-0}
   0 \le \int_{t_{j-1}}^{t_j} \left(t_j-s\right)\phi_1(T-s) \D s
    = \E^{-\gamma(T-t_j)} \rho_j,\quad
        \rho_j\coloneqq \left\{\kappa_1 \left[\tau_j - \left(T-t_j\right) \ln\left(1+\frac{\tau_j}{T-t_j}\right)
                          \right] + \kappa_1'\frac{\tau_j^2}{2}\right\} \,.
\end{gather}

Another example that appears frequently is, for $k=0,1,\dots$,
\begin{gather}\label{intG_t-1}
  \abs{\int_{t_{j-1}}^{t_j} \bigl(t_j-s\bigr)^k \bigl(s-t_{j-1}\bigr) \pt_t \m{G}(T-s) \D s}
     \le \E^{-\gamma(T-t_j)} \Phi_{k,j}\,, \\
  \intertext{where}\notag
    \Phi_{k,j} \coloneqq \kappa_1 \mu_{k,j} + \kappa_1' \frac{\tau_j^{k+2}}{(k+1)(k+2)}
    \quad\text{and}\quad
    \mu_{k,j} \coloneqq \int_{t_{j-1}}^{t_j} \frac{\bigl(t_j-s\bigr)^k \bigl(s-t_{j-1}\bigr)}{T-s} \D s\,.
\end{gather}
The last integral can be computed recursively: (with $\tau_j\coloneqq t_j-t_{j-1}$)
\begin{gather*}
   \mu_{0,j} = -\tau_j +(T-t_{j-1})\ln \left(1+\frac{\tau_j}{T-t_j}\right)\,, \ \ \
   \mu_{k,j} = \frac{\tau_j^{k+1}}{k(k+1)} + \left(t_j-T\right) \mu_{k-1,j}, \ \ k=1,2,\dots
\end{gather*}
However, when $t_j$ is close to $0$, destructive cancellation occures.
Then an alternative is to compute $\mu_{k,j}$ using a suitable truncation of the
series expansion
\begin{gather*}
   \mu_{k,j} 
     = \tau_j^{k+1} \sum_{\ell=1}^\infty \frac{(-1)^{\ell+1}}{(\ell+k)(\ell+k+1)}
          \left(\frac{\tau_j}{T-t_j}\right)^\ell.
\end{gather*}

Furthermore, for $k>0$
\begin{gather*}
  \int_{t_{j-1}}^{t_j} \bigl(t_j-s\bigr)^k \bigl(s-t_{j-1}\bigr) \pt_t \m{G}(T-s) \D s
    = \int_{t_{j-1}}^{t_j} \frac{\D}{\D s} \Bigl[\bigl(t_j-s\bigr)^k \bigl(s-t_{j-1}\bigr)\Bigr]
           \m{G}(T-s) \D s.
\end{gather*}
Application of~\eqref{source:ass}, gives the alternative bound
\begin{gather}\label{intG_t-2}
  \abs{\int_{t_{j-1}}^{t_j} \bigl(t_j-s\bigr)^k \bigl(s-t_{j-1}\bigr) \pt_t \m{G}(T-s) \D s}
    \le \E^{-\gamma(T-t_j)} \Phi_{k,j}^*, \quad
    \Phi_{k,j}^* \coloneqq
        \kappa_0
        \int_{t_{j-1}}^{t_j} \abs{\frac{\D}{\D s} \Bigl[\bigl(t_j-s\bigr)^k \bigl(s-t_{j-1}\bigr)\Bigr]}
         \D s\,.
\end{gather}
Combining~\eqref{intG_t-1} and~\eqref{intG_t-2} gives
\begin{gather}\label{intG_t-comb}
  \abs{\int_{t_{j-1}}^{t_j} \bigl(t_j-s\bigr)^k \bigl(s-t_{j-1}\bigr) \pt_t \m{G}(T-s) \D s}
    \le \E^{-\gamma(T-t_j)} \min\left\{\Phi_{k,j},\Phi_{k,j}^*\right\} \eqqcolon \Psi_{k,j}\,.
\end{gather}

\section{Test problem}
\label{sect:test}

Throughout the paper we shall give numerical results for the linear
reaction-diffusion equation
\begin{subequations}\label{testproblem}
\begin{alignat}{2}
 \pt_t u  - u_{xx} + (5x+6) u & = \E^{-4t} - \cos\bigl(\pi(x+t)^3\bigr)\,,
            &\quad& \text{in} \quad (-1,1) \times (0,1],\\
 \intertext{subject to the initial condition}
   u(x,0) & = u^0(x)=\sin \frac{\pi(1+x)}{2}\,, && \text{for} \quad x\in[-1,1], \\
 \intertext{and the Dirichlet boundary condition}
   u(x,t) & = 0\,, && \text{for} \quad (x,t)\in \{-1,1\} \times [0,1].
\end{alignat}
\end{subequations}
The Green's function for this problem satisfies~\cite[Corollary 5]{MR1423289}
\begin{gather*}
  \norm{\m{G}(t)}_{1,\Omega} \le \E^{-t/2},
    \quad
  \norm{\pt_t^p\m{G}(t)}_{1,\Omega}
            \le \frac{3}{2^{3/2}} \frac{p! 18^{p-1}}{t^p} \,\E^{-t/2}, \ \
        p \in \{1,2\}
\end{gather*}

The elliptic problems obtained after semi-discretisation in time are solved
using a spectral collocation method with polynomials of degree $31$. This
allows to solve those problems almost to machine accuracy.
We are interested in the errors and error estimates at final time $T$.
A reference solution is computed using dG(2) in time.
This is a method of order $5$, cf.~\cite{MR1227985,MR826227}.

\section{The implicit Euler method}
\label{sect:euler}

We consider the first-order backward Euler discretisations in time applied
to problem~\eqref{problem}.
Let an arbitrary mesh in time be given by
\begin{gather*}
   \omega_t \colon 0=t_0<t_1<\cdots<t_M=T.
\end{gather*}
For \mbox{$j=1,\dots,M$} we set
\begin{gather*}
   I_j\coloneqq\bigl(t_{j-1},t_j\bigr)\,,
 \quad
   \tau_j\coloneqq t_j-t_{j-1}
   \quad \text{and} \quad
   \tau \coloneqq \max_{j=1,\dots,M} \tau_j\,.
\end{gather*}
Furthermore, for $\varsigma\in[0,1]$ let
%\begin{gather*}
\mbox{$t_{j-\varsigma}\coloneqq t_j-\varsigma\tau_j$} \ and \
\mbox{$v^{j-\varsigma}\coloneqq v(t_{j-\varsigma})$}.
%\end{gather*}

We discretise the abstract parabolic problem \eqref{problem} in time
on the mesh $\omega_t$ using the first-order backward Euler method as follows.
We associate an approximate solution \mbox{$U^j\in H_0^1(\Omega)$}
with the time level $t_j$ and require it to satisfy
\begin{gather}\label{Euler_method}
  \delta_t U^j + \m{L} U^j =f^j
  \quad\text{in}   \; \Omega,\quad j=1,\dots, M; \quad U^0=u^0,
\end{gather}
where
\begin{gather*}
    \delta_t U^j \coloneqq \frac{U^j-U^{j-1}}{\tau_j}
              \quad\text{and}\quad
    f^j \coloneqq f(\cdot,t_j).
\end{gather*}

\paragraph{\S4.1\! }
The central idea is to extend $U^j$ to a piecewise linear function $\hat{U}$
that is defined on all of the interval $[0,T]$,
and then invoke~\eqref{green-rep} with $w=u-\hat{U}$.
To this end, for any function $v$ defined on $\omega_t$, $t_j\mapsto v^j$,
we denote by $\hat{v}$ its piecewise linear interpolant, i.e.,
\begin{gather*}\label{interpolant}
  \hat{v}(s) \coloneqq
      v^j - \bigl(t_j-s\bigr) \delta_t v^j
    = v^{j-1} + \bigl(s-t_{j-1}\bigr) \delta_t v^j
    = v^{j-1/2} + \bigl(s-t_{j-1/2}\bigr) \delta_t v^j\,,
  \quad s\in\hat{I}_j,\ \ j=1,\dots, M.
\end{gather*}
Note that, $\pt_t \hat{v}(s)=\delta_t v^j$ for $s\in I_j$,
$j=1,\dots,M$.
Hence,
\begin{gather*}
  \pt_t \hat{U} = \delta_t U^j = f^j-\m{L}U^j \quad \text{in} \ \ I_j\,,
   \ \ \text{by~\eqref{Euler_method}.}
\end{gather*}
Recalling~\eqref{problem}, the residual of $\hat{U}$ in the differential
equation admits the representation
\begin{gather}\label{euler:residuum1}
  \left(\m{K}\bigl(u-\hat{U}\bigr)\right)(s)
    = f(s) - \pt_t \hat{U}(s) - \m{L}\left(U^j - \bigl(t_j-s\bigr) \delta_t U^j \right)
     = f(s)-f^j + \bigl(t_j-s\bigr) \, \delta_t \left(\m{L}U\right)^j,
      \quad s\in I_j.
\end{gather}
Invoking~\eqref{green-rep}, we obtain for the error at final time $T=t_M$
\begin{align}\label{euler:error1}
  \notag
  & u(x,T)-U^M(x) = \bigl(u-\hat{U}\bigr)(x,T) \\
  & \qquad
     = \sum_{j=1}^M
         \left\{
           \int_{I_j} \scal{\m{G}(T-s)}{f(s)-f^j} \D s +
           \int_{I_j} \bigl(t_j-s\bigr)\, \scal{\m{G}(T-s)}{\delta_t \left(\m{L}U\right)^j} \D s
         \right\} \\
  & \qquad \label{KL1-repr}
     = \sum_{j=1}^M
         \left\{
           \int_{I_j} \scal{\m{G}(T-s)}{f(s)-f^j} \D s +
           \int_{I_j} \bigl(t_j-s\bigr)\, \scal{\pt_t\m{G}(T-s)}{\delta_t U^j} \D s
         \right\},
\end{align}
because \mbox{$\left(\pt_t + \m{L}^*\right)\m{G}=0$}.
Using the H\"older inequality and~\eqref{source:ass}, we obtain two bounds:
\begin{subequations}\label{DLM-KL1}
\begin{align}
  \norm{u(T)-U^M}_{\infty,\Omega}
    & \le \sum_{j=1}^M
           \left\{
             \int_{I_j} \phi_0(T-s) \norm{f(s)-f^j}_{\infty,\Omega} \D s +
             \int_{I_j} \bigl(t_j-s\bigr)\, \phi_0(T-s) \D s
                          \norm{\delta_t \left(\m{L}U\right)^j}_{\infty,\Omega}
           \right\} \\
  \intertext{and}
  \norm{u(T)-U^M}_{\infty,\Omega}
    & \le \sum_{j=1}^M
           \left\{
             \int_{I_j} \phi_0(T-s) \norm{f(s)-f^j}_{\infty,\Omega} \D s +
             \int_{I_j} \phi_1(T-s) \bigl(t_j-s\bigr) \D s \,
                        \norm{\delta_t U^j}_{\infty,\Omega}
           \right\}.
\end{align}
\end{subequations}
Upon noting that the $\phi_i$, \mbox{$i=0,1$}, are non-increasing,
we obtain the following theorems.
The first resembles the result given in~\cite[\S4.3, Theorem 4.2]{MR2519598},
while the second was derived in~\cite[\S4, Theorem 4.1]{MR3056758}.
A version of the latter is also given in~\cite[\S1, Theorem 1.3]{MR1335652},
but without providing a proof and without fixing the constants.
\begin{theorem}\label{theo:DLM}
  The maximum-norm error of the backward Euler time
  discretisation~\eqref{Euler_method} satisfies the a posteriori bound
  \begin{gather*}
    \norm{u(T)-U^M}_{\infty,\Omega} \le
      \sum_{j=1}^M \E^{-\gamma (T-t_j)}
                   \left( \eta_{\bar{f}}^j + \eta_{\delta \m{L}U}^j \right)
  \end{gather*}
  with
  \begin{align*}
    \eta_{\bar{f}}^j \coloneqq
       \kappa_0 \int_{I_j} \norm{f(s)-f^j}_{\infty,\Omega} \D s
    \quad\text{and}\quad
    \eta_{\delta\m{L}U}^j \coloneqq
       \frac{\kappa_0\tau_j^2}{2} \norm{\delta_t \left(\m{L}U\right)^j}_{\infty,\Omega}.
  \end{align*}
\end{theorem}
%the following
%theorem which resembles the result given in~\cite[\S4, Theorem 4.1]{MR3056758}.

%%%%% Kopteva/Linß

\begin{theorem}\label{theo:KL1}
  The maximum-norm error of the backward Euler time
  discretisation~\eqref{Euler_method} satisfies the a posteriori bound
  \begin{gather*}
    \norm{u(T)-U^M}_{\infty,\Omega} \le
      \sum_{j=1}^M \E^{-\gamma (T-t_j)}
                   \left( \eta_{\bar{f}}^j + \eta_{\delta U}^j \right)
  \end{gather*}
  with $\eta_{\bar{f}}^j$ as in Theorem~\ref{theo:DLM},
  $ \eta_{\delta U}^j \coloneqq
               \rho_j \norm{\delta_t U^j}_{\infty,\Omega}$
  and $\rho_j$ from~\eqref{intG_t-0}.
\end{theorem}

The derivation of Theorem~\ref{theo:KL1} in~\cite{MR3056758} uses a different,
global argument employing a piecewise constant and discontinuous interpolant
of the $U^j$.
In doing so, it passed unnoticed that these bounds can be combined by
locally taking,
for each $j=1,\dots,M$, the smaller of the two bounds
in~\eqref{DLM-KL1}.
We arrive at the following novel result.

\begin{theorem}\label{theo:DLMKL}
  The maximum-norm error of the backward Euler time
  discretisation~\eqref{Euler_method}
  satisfies the a posteriori bound
  \begin{gather*}
    \norm{u(T)-U^M}_{\infty,\Omega} \le
      \sum_{j=1}^M \E^{-\gamma (T-t_j)}
                   \left( \eta_{\bar{f}}^j + \eta^j_{\min} \right)\,,
     \quad\text{with} \quad
        \eta_{\min}^j \coloneqq
                   \min\left\{\eta_{\delta U}^j, \eta_{\delta \m{L} U}^j \right\}
  \end{gather*}
  and the notation from Theorems~\ref{theo:DLM} and~\ref{theo:KL1}.
\end{theorem}

\begin{remark}
  The integral defining $\eta_{\bar{f}}^j$ can (in general) not been evaluated
  exactly, but needs to be approximated.
  Possible options are
  \begin{align*}
    \int_{I_j} \norm{f(s)-f^j}_{\infty,\Omega} \D s
       & \approx \frac{\tau_j}{2} 
               \norm{f^{j-1}-f^j}_{\infty,\Omega} && \text{trapezium rule,} \\
    \int_{I_j} \norm{f(s)-f^j}_{\infty,\Omega} \D s
       & \approx \frac{\tau_j}{6} 
         \biggl\{
               \norm{f^{j-1}-f^j}_{\infty,\Omega}
               +4\norm{f^{j-1/2}-f^j}_{\infty,\Omega}
         \biggr\} && \text{Simpson's rule}
  \end{align*}
\end{remark}

\begin{table}
\centerline{%
\begin{tabular}{c|c|cc|cc|cc}
\multicolumn{2}{c}{} & \multicolumn{2}{|c}{Theorem 1} & \multicolumn{2}{|c}{Theorem 2} 
        & \multicolumn{2}{|c}{Theorem 3} \\\hline
$M$ & err & est & eff & est & eff & est & eff \\\hline
%    16 & 2.108e-03 & 1.069e+00 & 1/507 & 5.889e-01 & 1/279 & 5.889e-01 & 1/279 \\
%    32 & 9.395e-04 & 5.758e-01 & 1/613 & 3.081e-01 & 1/328 & 3.081e-01 & 1/328 \\
%    64 & 4.408e-04 & 2.931e-01 & 1/665 & 1.563e-01 & 1/355 & 1.563e-01 & 1/355 \\
%   128 & 2.129e-04 & 1.467e-01 & 1/689 & 7.877e-02 & 1/370 & 7.778e-02 & 1/365 \\
   256 & 1.045e-04 & 7.333e-02 & 1/702 & 3.970e-02 & 1/380 & 3.872e-02 & 1/370 \\
   512 & 5.175e-05 & 3.664e-02 & 1/708 & 2.003e-02 & 1/387 & 1.934e-02 & 1/374 \\
  1024 & 2.575e-05 & 1.831e-02 & 1/711 & 1.011e-02 & 1/393 & 9.662e-03 & 1/375 \\
  2048 & 1.284e-05 & 9.155e-03 & 1/713 & 5.106e-03 & 1/398 & 4.829e-03 & 1/376 \\
  4096 & 6.412e-06 & 4.577e-03 & 1/714 & 2.578e-03 & 1/402 & 2.414e-03 & 1/377 \\
  8192 & 3.204e-06 & 2.288e-03 & 1/714 & 1.302e-03 & 1/406 & 1.207e-03 & 1/377 \\
 16384 & 1.601e-06 & 1.144e-03 & 1/715 & 6.576e-04 & 1/411 & 6.035e-04 & 1/377 \\
 32768 & 8.006e-07 & 5.721e-04 & 1/715 & 3.320e-04 & 1/415 & 3.017e-04 & 1/377 \\
 65536 & 4.002e-07 & 2.860e-04 & 1/715 & 1.676e-04 & 1/419 & 1.509e-04 & 1/377 \\
\hline
\end{tabular}}
\caption{\label{tab:DLM-DLMKL}%
  Error estimators of Theorems~\ref{theo:DLM}-\ref{theo:DLMKL} applied to the test
  problem~\eqref{testproblem}.
  Simpson's rule is used to estimate the $\eta_{\bar{f}}^j$.
}
\end{table}

\paragraph{Numerical results.}
Table~\ref{tab:DLM-DLMKL} displays the results of our test computations
for~\eqref{testproblem}.
The first column contains the number of mesh intervals used on the spatial
domain $[0,1]$. To avoid special effects from uniform meshes, we have chosen
the mesh sizes to satisfy $\tau_j=2\tau_{j-1}$ for $j=2,4,6,\dots,M$.
The second column of the table displays the actual errors of the backward
Euler semidisretisation~\eqref{Euler_method}.
We observe convergence of order $1$ -- each time the number of mesh intervals
is doubled the error is divided by (approximately) two.

Columns 3 and 4 contain the a posteriori error bounds provided by
Theorem~\ref{theo:DLM} and its efficiency, i.e. the actual error divided by
the error estimator.
There is a strong correlation between the two.
However, the errors are overestimated by a factor of about $700$.

In columns 5 and 6 we have the corresponding numbers for Theorem~\ref{theo:KL1}.
It gives sharper bounds than Theorem~\ref{theo:DLM}, but the efficiency is slightly
deteriorating with the logarithm of the mesh size.
(Our test problem somewhat favours Theorem~\ref{theo:KL1}. There are other
equations where Theorem~\ref{theo:DLM} gives sharper bounds.)

Finally, in the last two columns of Table~\ref{tab:DLM-DLMKL} we present our
results for Theorem~\ref{theo:DLMKL}.
It gives sharper bounds than both Theorems~\ref{theo:DLM} and~\ref{theo:KL1},
which had to be expected from its derivation.
Moreover, we do not witness any deterioration of the efficiency with refinement
of the mesh.
Since the error bound of Theorem~\ref{theo:DLMKL} contains the mininum of two
terms, $\eta_{\delta U}^j$ and $\eta_{\delta \m{L} U}^j$, it is interesting to
study when which term is active. We will do this in a broader context later.

\paragraph{\S4.2\! }
The preceeding error bounds all contain a piecewise constant
approximation of the RHS $f$ of the PDE.
Now we shall involve its piecewise linear interpolation $\hat{f}$.
To this end we use $\hat{f}(t) = f^j - \bigl(t_j-t\bigr) \, \delta_t f^j$
and rewrite the residuum in~\eqref{euler:residuum1} as
\begin{gather*}
  \left(\m{K}\bigl(u-\hat{U}\bigr)\right)(t)
     = \bigl(f - \hat{f}\bigr)(t)
         + \bigl(t_j-t\bigr) \, \delta_t \left(\m{L}U-f\right)^j \ \ t\in I_j\,.
\end{gather*}
In view of~\eqref{Euler_method} we set \mbox{$\delta_t U^0 \coloneqq f^0-\m{L}U^0$},
introduce
\begin{gather*}
  \delta_t^2 v^j \coloneqq \frac{\delta_t v^j - \delta_t v^{j-1}}{\tau_j}\,,
    \quad j=1,\dots,M,
\end{gather*}
and obtain
\begin{align*}
  \left(\m{K}\bigl(u-\hat{U}\bigr)\right)(t)
     = \bigl(f - \hat{f}\bigr)(t)
         - \bigl(t_j-t\bigr) \, \delta_t^2 U^j\,, \quad t\in I_j, \ \ j=1,\dots,M.
\end{align*}
Proceeding as before, we get
\begin{theorem}\label{theo:LR}
  The maximum-norm error of the backward Euler time
  discretisation~\eqref{Euler_method} satisfies the a posteriori bound
  \begin{gather*}
    \norm{u(T)-U^M}_{\infty,\Omega}
      \le \sum_{j=1}^M \E^{-\gamma (T-t_j)}
                       \left(\eta_{\hat{f}}^j + \eta_{\delta^2 U}^j\right)
  \end{gather*}
  with
  \begin{align*}
    \eta_{\hat{f}}^j \coloneqq
       \kappa_0 \int_{I_j} \norm{\bigl(f-\hat{f}\bigr)(s)}_{\infty,\Omega} \D s\,, \quad
    \eta_{\delta^2 U}^j \coloneqq
       \frac{\kappa_0\tau_j^2}{2} \norm{\delta_t^2 U^j}_{\infty,\Omega} .
  \end{align*}
\end{theorem}
\begin{remark}
  Again, the integrals composing $\eta_{\hat{f}}$ need to be approximated.
  This time the trapizium rule would always give zero.
  One possibility is  Simpson's rule which gives
  \begin{gather*}
    \int_{I_j} \norm{\bigl(f-\hat{f}\bigr)(s)}_{\infty,\Omega} \D s
       \approx \frac{2\tau_j}{3} 
            \norm{\bigl(\hat{f}-f\bigr)^{j-1/2}}_{\infty,\Omega}
       = \frac{\tau_j}{3} 
            \norm{f^j - 2 f^{j-1/2} + f^{j-1}}_{\infty,\Omega}\,.
  \end{gather*}
\end{remark}

Taking minima locally for each time level $j$, $j=1,\dots,M$,
Theorems~\ref{theo:DLMKL} and Theorem~\ref{theo:LR} can be combined to give
the following sharpend result.
\begin{theorem}\label{theo:combined}
  The maximum-norm error of the backward Euler time
  discretisation~\eqref{Euler_method} satisfies the a posteriori bound
  \begin{gather*}
    \norm{u(T)-U^M}_{\infty,\Omega} \le
      \sum_{j=1}^M \E^{-\gamma (T-t_j)}
        \min\left\{ \eta_{\bar{f}}^j + \eta_{\min}^j,
                       \eta_{\hat{f}}^j + \eta_{\delta^2 U}^j\right\}
  \end{gather*}
  %
  %with $\eta_{\bar{f}}^j$ from Theorem~\ref{theo:DLM}, $\eta_{\min}^j$ from
  %Theorem~\ref{theo:DLMKL} and $\eta_{\hat{f}}^j + \eta_{\delta^2 U}^j$ from
  %Theorem~\ref{theo:LR}.
  with the notation from Theorems~\ref{theo:DLM}-\ref{theo:LR}.
\end{theorem}

\begin{table}
\centerline{%
\begin{tabular}{c|c|cc|cc}
\multicolumn{2}{c}{} & \multicolumn{2}{|c}{Theorem 4}
                     & \multicolumn{2}{|c}{Theorem 5} \\\hline
$M$ & err & est & eff & est & eff \\\hline
   256 & 1.045e-04 & 9.900e-03 & 1/95 & 6.498e-03 & 1/62 \\
   512 & 5.175e-05 & 4.796e-03 & 1/93 & 3.079e-03 & 1/59 \\
  1024 & 2.575e-05 & 2.360e-03 & 1/92 & 1.498e-03 & 1/58 \\
  2048 & 1.284e-05 & 1.171e-03 & 1/91 & 7.389e-04 & 1/58 \\
  4096 & 6.412e-06 & 5.834e-04 & 1/91 & 3.670e-04 & 1/57 \\
  8192 & 3.204e-06 & 2.911e-04 & 1/91 & 1.829e-04 & 1/57 \\
 16384 & 1.601e-06 & 1.454e-04 & 1/91 & 9.130e-05 & 1/57 \\
 32768 & 8.006e-07 & 7.269e-05 & 1/91 & 4.561e-05 & 1/57 \\
 65536 & 4.002e-07 & 3.634e-05 & 1/91 & 2.280e-05 & 1/57 \\
\hline
\end{tabular}}
\caption{\label{tab:combined}%
  Error estimators of Theorems~\ref{theo:LR} and~\ref{theo:combined}
  applied~\eqref{testproblem}.
  Simpson's rule is used to approximate $\eta_{\hat{f}}^j$,
  $\eta_{\bar{f}}^j$.}
\end{table}

\paragraph{Numerical results and discussions.}
Table~\ref{tab:combined} contains our results for Theorems~\ref{theo:LR}
and~\ref{theo:combined}. Both give sharper bounds than
Theorems~\ref{theo:DLM}-\ref{theo:DLMKL}.
This was expected for Theorem~\ref{theo:combined}.

\begin{figure}[b]
\centerline{\includegraphics[width=.6\textwidth]{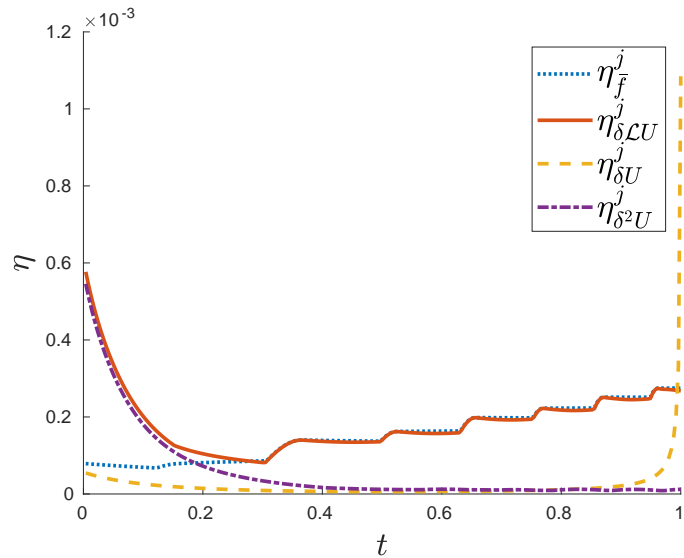}}
\caption{\label{fig:euler}
  The various parts of the error estimators in Theorems~\ref{theo:DLM}-%
  \ref{theo:combined}, uniform time stepping, $M=256$ steps.
}
\end{figure}

How do the various components of the error estimators behave?
Figure~\ref{fig:euler} depicts plots of the four terms
$\eta_{\bar{f}}\,$, $\eta_{\delta \m{L}U}$, $\eta_{\delta U}$ and
$\eta_{\delta^2 U}$.
We have chosen a uniform mesh as otherwise there would be oscillations because
the components are correlated with powers of the local mesh step size.
Also the term $\eta_{\hat{f}}\,$ is omitted because it is of higher order
and close to zero.
For the same reason graphs of
$\eta_{\delta^2 U}$ and $\eta_{\delta^2 U}+ \eta_{\hat{f}}\,$
would be virtually undistinguishable.

First, we notice that $\eta_{\delta \m{L}U}$ and $\eta_{\delta^2 U}$ attain
large values near inital time.
Second, $\eta_{\delta U}$ becomes large towards the final time. This can be
explained by the behaviour of the $\mu_j^\mathrm{bE}$ introduced in Theorem~\ref{theo:KL1}.
At final time \mbox{$t_M=T$}, we have \mbox{$\mu_M^\mathrm{bE}=\tau_M$}.
But further back in time, it becomes second order: \mbox{$\mu_j^\mathrm{bE}\sim\tau_j^2$}.

Theorems~\ref{theo:DLM} and~\ref{theo:KL1} differ in the use of
$\eta_{\delta\m{L}U}$ (solid red line) and $\eta_{\delta U}$ (dashed yellow
line). On most of the domain we have
\mbox{$\eta_{\delta\m{L}U}^j < \eta_{\delta U}^j$} only for the last few steps
the relation is reversed.
This illustrates how Theorem~\ref{theo:DLMKL} takes advantage by picking the
minimum of the two at each time step.

Finally, one notices that for times \mbox{$t\ge 0.3$} the terms $\eta_{\bar{f}}$
and $\eta_{\delta\m{L}U}$ take very similar values.
This suggests that in deriving Theorem~\ref{theo:DLM} a triangle inequality
might have been applied inadequately.
To illustrate this we look at the two representations of the residuum used above:
\begin{gather*}
  \underbrace{f(s)-f^j}_{\to \ \eta_{\bar{f}}^j} +
  \underbrace{\bigl(t_j-s\bigr)\,\delta_t \left(\m{L}U\right)^j}_{\to \ \eta_{\delta\m{L}U}^j}
    =
  \underbrace{\bigl(f-\hat{f}\bigr)(s)}_{\to \ \eta_{\hat{f}}^j} +
  \underbrace{\bigl(t_j-s\bigr)\,\delta_t\bigl(\m{L}U-f\bigr)^j}_{\to \ \eta_{\delta^2U}^j}\,.
\end{gather*}
Generically, the term $\eta_{\hat{f}}^j$ is of order $3$ (in $\tau_j$), while
the other three terms are of order $2$ only.
Therefore, asymptotically we have
\begin{gather*}
  \eta_{\delta^2U}^j \le \left(1+\m{O}\bigl(\tau_j\bigr)\right)
                         \left(\eta_{\bar{f}}^j + \eta_{\delta\m{L}U}^j\right)
                     \quad \bigl(\tau_j\to0\bigr)\,.
\end{gather*}
Thus, in general Theorem~\ref{theo:LR} will give sharper bounds than
Theorem~\ref{theo:DLM}.

In practice Theorem~\ref{theo:combined} should be given preference as it gives
the sharpest error bound.

\paragraph{\S4.3\! }
Concluding our study of the backward-Euler scheme, we like to review
an idea presented in~\cite{MR3720388}.
The primary intention of the authors was to eliminate the logarithmic
dependence on the time-step size observed in Theorem~\ref{theo:KL1}.

Let
\begin{gather*}
  W^j \coloneqq \frac{1}{2}
                \left[ \tau_j \delta_t U^j - \tau_M \delta_t U^M \right],
   \quad j=1,\dots,M.
\end{gather*}
The expectation in~\cite{MR3720388} was that for \mbox{$j\to M$} the $W^j$ behave
similar to \mbox{$T-t_j$}, and therefore compensate for the term \mbox{$T-s$} in the
denominator of the bound $\phi_1$ on $\m{G}_t$.
Then
\begin{gather*}
  \bigl(t_j - s\bigr)\, \delta_t U^j
    = \frac{\tau_M}{2} \delta_t U^M + W^j
         + \bigl(t_{j-1/2}-s\bigr)\, \delta_t U^j\,,
       \quad s\in\bigl(t_{j-1},t_j\bigr], \ \ j=1,\dots,M.
\end{gather*}
Define
\begin{gather*}
  \omega(s) \coloneqq \frac{\bigl(t_j-s\bigr)\bigl(s-t_{j-1}\bigr)}{2},
       \ \ s\in\bar{I}_j, \ \ j=1,\dots,M,
\end{gather*}
and note that
\begin{gather*}
  t_{j-1/2}-s = -\omega'(s) , \ \ s\in I_j\,.
\end{gather*}
Fix \mbox{$J\in\{1,\dots,M\}$}.
Integration by parts for the interval \mbox{$[t_{J-1},t_{M-1}]$} applied to
the second term on the RHS of~\eqref{KL1-repr} gives
\begin{align*}
  u(x,T)-U^M(x)
   &  = \sum_{j=1}^M \int_{I_j} \scal{\m{G}(T-s)}{f(s)-f^j} \D s
        + \sum_{j=1}^{J-1,M} \int_{I_j} \bigl(t_j-s\bigr)\, \scal{\pt_t\m{G}(T-s)}{\delta_t U^j} \D s \\
   & \qquad
        - \sum_{j=J}^{M-1}
            \left\{\int_{I_j} \omega(s) \scal{\pt_t^2 \m{G}(T-s)}{\delta_t U^{j}} \D s
                 - \int_{I_j} \scal{\pt_t \m{G}(T-s)}{W^j} \D s
            \right\} \\
   & \qquad
        - \frac{\tau_M}{2} \scal{\m{G}(T-t_{M-1}) - \m{G}(T-t_{J-1})}{\delta_t U^{M}}.
\end{align*}
The first and second integral are estimated as in the derivation of Theorem~\ref{theo:DLMKL}.
To the third and fourth integral we apply~\eqref{intG_p}.
The last one is bounded using H\"older's inequality again and~\eqref{source:ass}.
%We arrive at
%%
%\begin{align*}
%  \norm{u(T)-U^M}_{\infty,\Omega}
%   & \le \sum_{j=1}^M \E^{-\gamma(T-t_j)} \eta_{\bar{f}}^j
%           + \sum_{j=1}^{J-1,M} \E^{-\gamma(T-t_j)}
%                                \min\left\{\eta_{\delta U}^j,\eta_{\delta\m{L}U}^j
%                                    \right\} \\
%    & \qquad
%       + \sum_{j=J}^{M-1}
%           \left\{   \int_{I_j} \omega(s) \phi_2(T-s) \D s \,
%                     \norm{\delta_t U^j}_{\infty,\Omega}
%                   + \int_{I_j} \phi_1(T-s) \D s \, \norm{W^j}_{\infty,\Omega}
%           \right\} \\
%   & \qquad
%        + \frac{\tau_M}{2} \Bigl(\phi_0(T-t_{M-1}) + \phi_0(T-t_{J-1})\Bigr)
%               \norm{\delta_t U^M}_{\infty,\Omega}\,.
%\end{align*}
%
\begin{theorem}\label{theo:KL2}
  For any $J\in\{1,\dots,M\}$ the maximum-norm error of the backward
  Euler time discretisation~\eqref{Euler_method}
  satisfies the a posteriori bound
  \begin{align*}
    \norm{u(T)-U^M}_{\infty,\Omega} & \le
      \sum_{j=1}^M \E^{-\gamma (T-t_j)} \eta_{\bar{f}}^j +
      \sum_{j=1}^{J-1,M} \E^{-\gamma (T-t_j)}
                         \min\left\{\eta_{\delta U}^j,\eta_{\delta\m{L}U}^j
                             \right\} \\
    & \qquad
       + \sum_{j=J}^{M-1} \E^{-\gamma (T-t_j)}
           \left\{ \eta_{\delta U,*}^j + \eta_W^j
           \right\}
        + \frac{\kappa_0\tau_M}{2} \Bigl(\E^{-\gamma(T-t_{M-1})} + \E^{-\gamma(T-t_{J-1})}\Bigr)
               \norm{\delta_t U^M}_{\infty,\Omega}\,.
  \end{align*}
  with $\eta_{\bar{f}}^j$ and $\eta_{\delta\m{L}U}^j$ from Theorem~\ref{theo:DLM}
  and $\eta_{\delta U}^j$ from Theorem~\ref{theo:KL1} and the new terms
  \begin{gather*}
    \eta_{\delta U,*}^j \coloneqq
       \left(\kappa_2 \mu_j^* + \frac{\kappa_2' \tau_j^3}{6}\right)
               \norm{\delta_t U^j}_{\infty,\Omega}\, , \quad
    \eta_W^j \coloneqq \theta_j % \left(\kappa_1 \ln \frac{T-t_{j-1}}{T-t_j} + \kappa_1'\tau_j\right)
                \norm{W^j}_{\infty,\Omega}\,,\quad
    \mu_j^* \coloneqq \int_{I_j} \frac{\omega(s)}{(T-s)^2} \D s\,.
  \end{gather*}
\end{theorem}
\begin{remark}\label{rem:bE:doublerun}
  In~\cite{MR3720388} the result is derived for \mbox{$J=1$} and with only
  $\eta_{\delta U}^j$ in the second sum instead of
  \mbox{$\min\left\{\eta_{\delta U}^j,\eta_{\delta\m{L}U}^j\right\}$}.

  The drawback of this approach is that in order to compute the $W^j$ one has
  to know $U^M$ and $U^{M-1}$. Hence, one either has to perform two runs for
  \mbox{$j=J,\dots,M$}, the first to determine $\delta_tU^M$ and the second
  to compute the $W^j$, or one
  needs to store the approximations at those time levels.
\end{remark}

\begin{table}
\centerline{%
\begin{tabular}{c|c|cc}
\multicolumn{2}{c}{} & \multicolumn{2}{|c}{Theorem 6} \\\hline
$M$ & err & est & eff \\\hline
%    16 & 2.108e-03 & 5.825e-01 & 1/276 \\
%    32 & 9.395e-04 & 2.991e-01 & 1/318 \\
%    64 & 4.408e-04 & 1.494e-01 & 1/339 \\
%   128 & 2.129e-04 & 7.326e-02 & 1/344 \\
   256 & 1.045e-04 & 3.596e-02 & 1/344 \\
   512 & 5.175e-05 & 1.778e-02 & 1/344 \\
  1024 & 2.575e-05 & 8.833e-03 & 1/343 \\
  2048 & 1.284e-05 & 4.396e-03 & 1/342 \\
  4096 & 6.412e-06 & 2.190e-03 & 1/342 \\
  8192 & 3.204e-06 & 1.092e-03 & 1/341 \\
 16384 & 1.601e-06 & 5.446e-04 & 1/340 \\
 32768 & 8.006e-07 & 2.716e-04 & 1/339 \\
 65536 & 4.002e-07 & 1.355e-04 & 1/338 \\
\hline
\end{tabular}}
\caption{\label{tab:KL2}%
  Error estimator of Theorem~\ref{theo:KL2}, \mbox{$J=1$} applied to the test
  problem~\eqref{testproblem}.
  Simpson's rule is used again to estimate the $\eta_{\bar{f}}^j$.
}
\end{table}

\paragraph{Numerical results.}
Table~\ref{tab:KL2} displays our numerical results for Theorem~\ref{theo:KL2}.
We witness a slight improvement over the error bounds of Theorems~\ref{theo:DLMKL},
but not over Theorem~\ref{theo:combined}.

\section{The Crank-Nicolson method}
\label{sect:CN}

We discretise the abstract parabolic problem \eqref{problem} in time
on the mesh $\omega_t$ using the second-order Crank-Nicolson method as follows.
We associate an approximate solution \mbox{$U^j\in H_0^1(\Omega)$}
with the time level $t_j$ and require it to satisfy
\begin{align}\label{CN_method}
  \delta_t U^j + \m{L} \hat{U}^{j-1/2} & = \hat{f}^{j-1/2}
  \quad\text{in}   \; \Omega,\quad j=1,\dots, M; \quad U^0=u^0\,,
  \intertext{i.e.} \notag
  \frac{U^j-U^{j-1}}{\tau_j} + \frac{\m{L} U^j + \m{L} U^{j-1}}{2}
     & = \frac{f^j + f^{j-1}}{2}
  \quad\text{in}   \; \Omega,\quad j=1,\dots, M; \quad U^0=u^0\,.
\end{align}

\paragraph{\S5.1\! }
We extend the $U^j$ to a globally defined function using piecewise linear
interpolation:
\begin{gather*}
  \hat{U}(s) = U^j - \bigl(t_j-s\bigr) \, \delta_t U^j
    = \hat{U}^{j-1/2} + \bigl(s-t_{j-1/2}\bigr) \, \delta_t U^j\,,
  \quad s\in \hat{I}_j \,,\ \ j=1,\dots, M.
\end{gather*}
The residuum of $\hat{U}$ in the PDE admits the representation
\begin{gather*}
  \left(\m{K}\bigl(u-\hat{U}\big)\right)(s)
    = f(s) -\pt_t \hat{U}(s)
           - \m{L} \left(\hat{U}^{j-1/2} + \bigl(s-t_{j-1/2}\bigr)\,
                   \delta_t U^j\right),\quad s\in I_j.
\end{gather*}
Let \mbox{$\psi^j \coloneqq \bigl(\m{L}U-f\bigr)^j$}.
Then by~\eqref{CN_method}, we have \
\mbox{$\pt_t \hat{U}(s) = \delta_t U^j
        = \hat{f}^{j-1/2} - \m{L}\hat{U}^{j-1/2}
        = - \hat{\psi}^{j-1/2}$}
for \mbox{$s\in I_j$}.
This gives
\begin{gather*}
  \left(\m{K}\bigl(u-\hat{U}\bigr)\right)(s)
    = f(s) - \hat{f}^{j-1/2}
           + \bigl(t_{j-1/2}-s\bigr) \, \delta_t\bigl(\m{L} U\bigr)^j
    = f(s) - \hat{f}(s)
           + \bigl(t_{j-1/2}-s\bigr) \, \delta_t\psi^j\,,\quad s\in I_j.
\end{gather*}
We substitute into~\eqref{green-rep} and obtain
\begin{gather}\label{CN-repr}
  u(x,T) - U^M(x)
    = \sum_{j=1}^{M}
      \left\{   \int_{I_j} \scal{\m{G}(T-s)}{\bigl(f-\hat{f}\bigr)(s)}\D s
              + \int_{I_j} \bigl(t_{j-1/2}-s\bigr)
                \scal{\m{G}(T-s)}{\delta_t\psi^j}\D s
      \right\}.
\end{gather}
To the first integral we apply~\eqref{intG}.
When bounding the second one, note that
$\bigl(t_{j-1/2}-s\bigr) = \frac{1}{2}\frac{\D}{\D s} \bigl(t_j-s\bigr)\bigl(s-t_{j-1}\bigr)$.
Therefore, we can avail of~\eqref{intG_t-comb} for $k=1$.
%
%Using~\eqref{source:ass}, we get a first bound for the terms on the RHS:
%%
%\begin{gather*}
%  \abs{\int_{I_j} \scal{\m{G}(T-s)}{\bigl(f-\hat{f}\bigr)(s)}\D s}
%      \le \kappa_0 \E^{-\gamma(T-t_j)}
%          \int_{I_j} \norm{\bigl(f-\hat{f}\bigr)(s)}_{\infty,\Omega} \D s
%\end{gather*}
%%
%and
%%
%\begin{gather}\label{CN-dpsi-0}
%  \abs{\int_{I_j} \bigl(t_{j-1/2}-s\bigr)
%                \scal{\m{G}(T-s)}{\delta_t\psi^j} \D s}
%      \le \frac{\kappa_0 \tau_j^2}{4} \E^{-\gamma(T-t_j)}
%          \norm{\delta_t\psi^j}_{\infty,\Omega}
%\end{gather}
%%
%An alternative bound is obtained using integration by parts
%-- recall that $\omega(s)=\frac{1}{2}\bigl(t_j-s\bigr)\bigl(s-t_{j-1}\bigr)$ --
%%
%\begin{gather}\label{CN-repr-*}
%  \int_{I_j} \bigl(t_{j-1/2}-s\bigr) \scal{\m{G}(T-s)}{\delta_t\psi^j} \D s
%    = \int_{I_j} \omega(s) \scal{\pt_t\m{G}(T-s)}{\delta_t\psi^j}\D s
%\end{gather}
%%
%and estimating as follows
%%
%\begin{gather*}
%  \abs{\int_{I_j} \bigl(t_{j-1/2}-s\bigr)
%                \scal{\m{G}(T-s)}{\delta_t\psi^j} \D s}
%      \le \int_{I_j} \omega(s)
%                   \left(\frac{\kappa_1}{T-s}+\kappa_1'\right) \D s \ 
%                   \E^{-\gamma(T-t_j)} \norm{\delta_t\psi^j}_{\infty,\Omega}\,,
%\end{gather*}
%%
%by~\eqref{source:ass}, $p=1$.
We arrive at the following theorem which is a slight modification of the result given
in~\cite[\S5, Theorem 5.1]{MR3056758}.
\begin{theorem}\label{theo:CN1}
  The maximum-norm error of the Crank-Nicolson method~\eqref{CN_method}
  satisfies the a posteriori error bound
  \begin{gather*}
    \norm{u(T)-U^M}_{\infty,\Omega}
      \le \sum_{j=1}^M \E^{-\gamma (T-t_j)}
                       \left(\eta_{\hat{f}}^j + \eta_{\delta\psi}^j\right)
  \end{gather*}
  with $\eta_{\hat{f}}^j$ as in Theorem~\ref{theo:LR},
  \begin{gather*}
    \eta_{\delta\psi}^j \coloneqq
         \frac{\Psi_{1,j}}{2} \norm{\delta_t\psi^j}_{\infty,\Omega}\,,\quad
    \psi^j \coloneqq \bigl(\m{L}U-f\bigr)^j
  \end{gather*}
  and $\Psi_{1,j}$ from~\eqref{intG_t-comb}.
\end{theorem}

\paragraph{\S5.2\! }
When studying the backward Euler semidiscretisation, the use of a
higher order interpolant of the RHS $f$ turned out to be useful.
This time, we define a piecewise quadratic interpolant $\tilde{f}$ by
\begin{gather*}
  \tilde{f}(s) \coloneqq \hat{f}(s) + \beta_j \omega(s), \quad s\in \bar{I_j},
     \quad\text{with}\quad
     \beta_j \coloneqq - 4 \frac{f^j - 2f^{j-1/2} + f^{j-1}}{\tau_j^2}
        \approx - \bigl(f''\bigr)^{j-1/2}\,.
\end{gather*}
It interpolates $f$ at the mesh points of $\omega_t$ and at the midpoint of its
mesh intervals.
Let \mbox{$\m{L}^{-1}\beta^j\coloneqq q^j\in H_0^1(\Omega)$} be
the unique solution of \mbox{$\m{L}q^j=\beta^j$}.
Then
\begin{align*}
  & \int_{I_j} \omega(s) \scal{\m{G}(T-s)}{\beta^j}\D s
    = \int_{I_j} \omega(s) \scal{\m{G}(T-s)}{\m{L}q^j}\D s \\
  & \qquad
    = - \int_{I_j} \omega(s) \scal{\pt_t\m{G}(T-s)}{q^j}\D s
    = - \int_{I_j} \omega'(s) \scal{\m{G}(T-s)}{q^j}\D s\,,
\end{align*}
because \mbox{$\m{L}^*\m{G}=-\pt_t\m{G}$}, and by integration by parts.
Then, from~\eqref{CN-repr}
\begin{align*}
  u(x,T) - U^M(x)
  & = \sum_{j=1}^{M}
      \left\{   \int_{I_j} \scal{\m{G}(T-s)}{\bigl(f-\tilde{f}\bigr)(s)}\D s
              + \int_{I_j} \omega'(s) \scal{\m{G}(T-s)}{\delta_t\psi^j-q^j}\D s
      \right\}.
% \\
%  & = \sum_{j=1}^{M}
%      \left\{   \int_{I_j} \scal{\m{G}(T-s)}{\bigl(f-\tilde{f}\bigr)(s)}\D s
%              + \int_{I_j} \omega(s) \scal{\pt_t\m{G}(T-s)}{\delta_t\psi^j-q^j}\D s
%      \right\}.
\end{align*}
Using the H\"older inequality, \eqref{intG} and~\eqref{intG_t-comb},
we obtain our next result.
\begin{theorem}\label{theo:CN2}
  The maximum-norm error of the Crank-Nicolson method~\eqref{CN_method}
  satisfies the a posteriori bound
  \begin{gather*}
    \norm{u(T)-U^M}_{\infty,\Omega}
      \le \sum_{j=1}^M \E^{-\gamma (T-t_j)}
                       \left(\eta_{\tilde{f}}^j + \eta_{\delta\psi q}^j\right)
  \end{gather*}
  with $q^j\in H_0^1(\Omega)$ solving \mbox{$\m{L}q^j=\beta^j$},
  \begin{gather*}
    \eta_{\tilde{f}}^j \coloneqq
       \kappa_0 \int_{I_j} \norm{\bigl(f-\tilde{f}\bigr)(s)}_{\infty,\Omega} \D s\,,\quad
    \eta_\chi^j \coloneqq
         \frac{\Psi_{1,j}}{2} \norm{\delta_t\psi^j-q^j}_{\infty,\Omega}\,
  \end{gather*}
  \mbox{$\psi^j\coloneqq \bigl(\m{L}U-f\bigr)^j$}
  and $\Psi_{1,j}$ from~\eqref{intG_t-comb}.
\end{theorem}
\begin{remark}
  The integral defining $\eta_{\tilde{f}}^j$ can (in general) not be evaluated
  exactly, but needs to be approximated.
  For example, Simpson's rule can be applied on the two subintervals
  \mbox{$[t_{j-1},t_{j-1/2}]$} and \mbox{$[t_{j-1/2},t_j]$} to give
  \begin{align*}
    \int_{I_j} \norm{\bigl(f-\tilde{f}\bigr)(s)}_{\infty,\Omega} \D s
       \approx \frac{\tau_j}{3} 
        \left\{ \norm{\bigl(f-\tilde{f}\bigr)^{j-3/4}}_{\infty,\Omega}
               +\norm{\bigl(f-\tilde{f}\bigr)^{j-1/4}}_{\infty,\Omega}
        \right\}.
  \end{align*}
\end{remark}

\begin{remark}
  The above choice of a piecewise quadratic interpolation of $f$ corresponds
  to a piecewise quadratic reconstruction
  \mbox{$\tilde{U}(s) = \hat{U}(s) + \m{L}^{-1}\beta^j \omega(s)$},
  of the approximations $U^j$.

  In~\cite{MR2196979} the authors also used a special piecewise quadratic
  reconstruction of the $U^j$ in an a posteriori error analysis, but in the
  context of error estimation in $L_2$-type norms.
\end{remark}

Again, taking minima locally for each time level $j$, \mbox{$j=1,\dots,M$},
the bounds of the previous two theorems can be combined to give the sharpened
result:
\begin{theorem}\label{theo:CNcomb}
  The maximum-norm error of the Crank-Nicolson method~\eqref{CN_method}
  satisfies the a posteriori bound
  \begin{gather*}
    \norm{u(T)-U^M}_{\infty,\Omega}
      \le \sum_{j=1}^M \E^{-\gamma (T-t_j)}
                       \min\left\{\eta_{\hat{f}}^j + \eta_{\delta\psi}^j,
                                  \eta_{\tilde{f}}^j + \eta_{\delta\psi q}^j
                           \right\},
  \end{gather*}
  with $\eta_{\hat{f}}^j$ as in Theorem~\ref{theo:LR},
  $\eta_{\delta\psi}^j$ in Theorem~\ref{theo:CN1} and
  $\eta_{\tilde{f}}^j$ and $\eta_{\delta\psi q}^j$ from Theorem~\ref{theo:CN2}.
\end{theorem}

\paragraph{Numerical results} for the Crank-Nicolson method are given in
Table~\ref{tab:CN}.
For our test problem, the estimator of Theorem~\ref{theo:CN1} overestimates the
errors by a factor of almost $2000$.
In contrast, Theorems~\ref{theo:CN2} and~\ref{theo:CNcomb} yield sharper error
bounds.
Of course with Theorem~\ref{theo:CNcomb} giving the best.
However, for all three the efficiency slightly deteriorates as the mehs is
refined.

\begin{table}
\centerline{%
\begin{tabular}{c|c|cc|cc|cc}
\multicolumn{2}{c}{} & \multicolumn{2}{|c}{Theorem \ref{theo:CN1}}
                     & \multicolumn{2}{|c}{Theorem \ref{theo:CN2}}
                     & \multicolumn{2}{|c}{Theorem \ref{theo:CNcomb}} \\\hline
$M$ & err & est & eff & est & eff & est & eff \\\hline
%    16 & 1.799e-04 & 1.782e-01 & 1/991 & 2.438e-01 & 1/1355 & 1.686e-01 & 1/937 \\
%    32 & 4.081e-05 & 4.559e-02 & 1/1117 & 3.256e-02 & 1/798 & 3.183e-02 & 1/780 \\
%    64 & 8.201e-06 & 1.150e-02 & 1/1402 & 4.466e-03 & 1/545 & 4.466e-03 & 1/545 \\
%   128 & 1.819e-06 & 2.891e-03 & 1/1589 & 6.583e-04 & 1/362 & 6.583e-04 & 1/362 \\
   256 & 4.284e-07 & 7.221e-04 & 1/1686 & 1.102e-04 & 1/257 & 1.102e-04 & 1/257 \\
   512 & 1.048e-07 & 1.806e-04 & 1/1724 & 2.134e-05 & 1/204 & 2.033e-05 & 1/194 \\
  1024 & 2.592e-08 & 4.529e-05 & 1/1747 & 4.744e-06 & 1/183 & 4.253e-06 & 1/164 \\
  2048 & 6.445e-09 & 1.145e-05 & 1/1776 & 1.245e-06 & 1/193 & 1.087e-06 & 1/169 \\
  4096 & 1.607e-09 & 2.875e-06 & 1/1789 & 3.239e-07 & 1/202 & 2.750e-07 & 1/171 \\
  8192 & 4.011e-10 & 7.221e-07 & 1/1800 & 8.631e-08 & 1/215 & 7.192e-08 & 1/179 \\
 16384 & 1.002e-10 & 1.814e-07 & 1/1810 & 2.318e-08 & 1/231 & 1.904e-08 & 1/190 \\
 32768 & 2.503e-11 & 4.557e-08 & 1/1821 & 6.222e-09 & 1/249 & 5.053e-09 & 1/202 \\
 65536 & 6.224e-12 & 1.145e-08 & 1/1839 & 1.666e-09 & 1/268 & 1.340e-09 & 1/215 \\
\hline
\end{tabular}}
\caption{\label{tab:CN}
  Error estimators of Theorems~\ref{theo:CN1}-\ref{theo:CNcomb}
  applied to test problem~\eqref{testproblem};
  Simpson's rule used to estimate $\eta_{\hat{f}}$.}
\end{table}

\paragraph{\S5.3\! }
Concluding our study of the Crank-Nicolson method, we review an idea
presented in~\cite{MR3720388}.

Let
\begin{gather*}
  W_\psi^j \coloneqq \frac{1}{12}
                \left[ \tau_j^2 \delta_t \psi^j - \tau_M^2 \delta_t \psi^M \right]
   \quad\text{and}\quad
   \tilde{\omega}(s) \coloneqq \omega(s) - \frac{\tau_j^2}{12}, \ \ s\in I_j,
   \quad j=1,\dots,M.
\end{gather*}
The expectation in~\cite{MR3720388} was that for \mbox{$j\to M$} the $W_\psi^j$ behave
similar to \mbox{$T-t_j$}, and therefore compensate for the term \mbox{$T-s$} in the
denominator of the bound $\phi_1$ on $\m{G}_t$, see~\eqref{source:ass}.
Then
\begin{gather*}
  \omega(s) \delta_t \psi^j
    = \frac{\tau_M^2}{12} \delta_t \psi^M + W_\psi^j
         + \tilde{\omega}(s) \, \delta_t \psi^j\,,
       \quad s\in \bar{I}_j, \ \ j=1,\dots,M.
\end{gather*}
Define
\begin{gather*}
  \pi(s) \coloneqq \int_{t_{j-1}}^s \tilde{\omega}(\sigma) \, \D \sigma
     = \frac{1}{6} \bigl(t_j-s\bigr)\bigl(t_{j-1/2}-s\bigr) \bigl(t_{j-1}-s\bigr)\,,
       \ \ s\in\bar{I}_j, \ \ j=1,\dots,M.
\end{gather*}
Fix \mbox{$J\in\{1,\dots,M\}$}.
Integration by parts applied to parts of the RHS of~\eqref{CN-repr} gives
\begin{align*}
  u(x,T)-U^M(x)
   &  = \sum_{j=1}^M \int_{I_j} \scal{\m{G}(T-s)}{\bigl(f-\hat{f}\bigr)(s)} \D s
        + \sum_{j=1}^{J-1,M} \int_{I_j} \omega(s) \scal{\pt_t\m{G}(T-s)}{\delta_t \psi^j} \D s \\
   & \qquad
        + \sum_{j=J}^{M-1}
            \left\{\int_{I_j} \pi(s) \scal{\pt_t^2 \m{G}(T-s)}{\delta_t \psi^j} \D s
                 + \int_{I_j} \scal{\pt_t \m{G}(T-s)}{W_\psi^j} \D s
            \right\} \\
   & \qquad
        + \frac{\tau_M^2}{12} \scal{\m{G}(T-t_{M-1}) - \m{G}(T-t_{J-1})}{\delta_t \psi^{M}}.
\end{align*}
We employ our standard machinery and arrive at
%Using H\"older's inequality and~\eqref{source:ass}, we arrive at
%%
%\begin{align*}
%  \norm{u(T)-U^M}_{\infty,\Omega}
%   & \le \sum_{j=1}^M \E^{-\gamma(T-t_j)} \eta_{\hat{f}}^j
%           + \sum_{j=1}^{J-1,M} \E^{-\gamma(T-t_j)} \eta_{\delta\psi}^j \\
%    & \qquad
%       + \sum_{j=J}^{M-1}
%           \left\{   \int_{I_j} \pi(s) \phi_2(T-s) \D s \,
%                     \norm{\delta_t \psi^j}_{\infty,\Omega}
%                   + \int_{I_j} \phi_1(T-s) \D s \, \norm{W_\psi^j}_{\infty,\Omega}
%           \right\} \\
%   & \qquad
%        + \frac{\tau_M}{2} \Bigl(\phi_0(T-t_{M-1}) + \phi_0(T-t_{J-1})\Bigr)
%               \norm{\delta_t \psi^M}_{\infty,\Omega}\,.
%\end{align*}
%
\begin{theorem}\label{theo:CN-KL2}
  For any $J\in\{1,\dots,M\}$ the maximum-norm error
  of the Crank-Nicolson method~\eqref{CN_method}
  satisfies the a posteriori bound
  \begin{align*}
    \norm{u(T)-U^M}_{\infty,\Omega} & \le
      \sum_{j=1}^M \E^{-\gamma (T-t_j)} \eta_{\hat{f}}^j +
      \sum_{j=1}^{J-1,M} \E^{-\gamma(T-t_j)} \eta_{\delta\psi}^j \\
    & \qquad
       + \sum_{j=J}^{M-1} \E^{-\gamma (T-t_j)}
           \left\{ \eta_{\delta\psi,*}^j + \eta_{W_\psi}^j \right\}
        + \frac{\kappa_0\tau_M^2}{12} \Bigl(\E^{-\gamma(T-t_{M-1})} + \E^{-\gamma(T-t_{J-1})}\Bigr)
               \norm{\delta_t \psi^M}_{\infty,\Omega}
  \end{align*}
  with $\eta_{\hat{f}}^j$ and $\eta_{\delta\psi}^j$ from Theorems~\ref{theo:LR} and~\ref{theo:CN1}
  and the new terms
  \begin{gather*}
    \eta_{\delta\psi,*}^j \coloneqq
       \left(\kappa_2 \sigma_j^* + \frac{\kappa_2' \tau_j^4}{144}\right)
       \norm{\delta_t \psi^j}_{\infty,\Omega}, \quad
    \eta_{W_\psi}^j \coloneqq \theta_j % \left(\kappa_1 \ln \frac{T-t_{j-1}}{T-t_j} + \kappa_1'\tau_j\right)
                \norm{W_\psi^j}_{\infty,\Omega},
    \quad
    \sigma_j^* \coloneqq \int_{I_j} \frac{\abs{\pi(s)}}{(T-s)^2} \D s .
  \end{gather*}
\end{theorem}
%
%\begin{remark}
  %In~\cite{MR3720388} the result is derived for \mbox{$J=1$}.
  %The drawback of this approach is that in order to compute the $W_\psi^j$ one has
  %to know $U^M$ and $U^{M-1}$. Hence, one either has to perform two complete runs,
  %the first to determine $\delta_tU^M$ and the second to compute the $W^j$, or one
  %needs to store the approximate solutions at all time levels.
%\end{remark}
Note, that Remark~\ref{rem:bE:doublerun} holds accordingly.

%
% Results for Theorem~\ref{theo:CN-KL2} are missing.

\section{Extrapolated Euler method}
\label{sect:extra}

This extrapolation method combines two approximations by the backward
Euler-method on the mesh $\omega_t$ and on a mesh that is twice as fine.
They are defined by
\begin{subequations}\label{euler-extra}
\begin{description*}
  \item[One-step Euler:] $V^0=u^0$,
     \begin{gather}\label{euler-one}
       \delta_t V^j + \m{L} V^j = f^j, \quad j=1,2,\dots,M,
     \end{gather}
  \item[Two-step Euler:] $W^0=u^0$,
     \begin{gather}\label{euler-two}
        \frac{W^{j-1/2}-W^{j-1}}{\tau_j/2} + \m{L}W^{j-1/2} = f^{j-1/2}, \quad
        \frac{W^j-W^{j-1/2}}{\tau_j/2} + \m{L}W^j = f^j, \quad j=1,\dots,M.
     \end{gather}
  \item[Extrapolation:]
     \begin{gather}\label{extra}
       U^j \coloneqq 2W^j - V^j, \quad j=1,\dots,M.
     \end{gather}
\end{description*}
\end{subequations}

We follow~\cite{RadLin22} and consider a piecewise linear reconstruction
$\hat{U}$ of the approximations $U^j$, \mbox{$j=0,1,\dots,M$}.
First, adding the two equations in \eqref{euler-two} and
subtracting~\eqref{euler-one} yields
\begin{gather*}
  \pt_t \hat{U} = \delta_t U^j = 2 \delta_t W^j - \delta_t V^j
         = f^{j-1/2} - \m{L}\left(W^{j-1/2} + W^j - V^j\right)\,.
\end{gather*}
This implies for the residuum
\begin{gather}\label{res-1}
   \left(\m{K}(u-\hat{U})\right)(s)
       = f(s) - \pt_t \hat{U}(s) - \m{L} \hat{U}(s)
       = f(s) - f^{j-1/2} + \m{L} \left(W^{j-1/2} - W^j\right) + \m{L}\left(U^j - \hat{U}(s)\right)\,.
\end{gather}
Next,
\begin{gather*}
  U^j - \hat{U}(s) = - (s-t_j) \delta_t U^j = -(s-t_{j-1/2}) \delta_t U^j
      + \frac{\tau_j}{2} \delta_t U^j\,,
\end{gather*}
which implies
\begin{gather*}
  \m{L}\left(U^j - \hat{U}(s)\right) = - (s-t_{j-1/2}) \m{L} \delta_t U^j
      + \frac{1}{2} \m{L} \left(U^j - U^{j-1}\right)\,.
\end{gather*}
This is substituted into~\eqref{res-1} to give
\begin{align*}
  \left(\m{K}(u-\hat{U})\right)(s)
     & = \bigl(f - \hat{f}\bigr)(s) + \hat{f}^{j-1/2} - f^{j-1/2}
             + (t_{j-1/2}-s) \, \delta_t \bigl(\m{L}U-f\bigr)^j \\
     & \qquad\quad
             +  \m{L}\left(W^{j-1/2} - W^{j-1} - \frac{V^j - V^{j-1}}{2}\right)\,.
\end{align*}
Setting
\begin{gather}\label{def:Zj}
   Z^j \coloneqq  W^{j-1/2} - W^{j-1} - \frac{V^j - V^{j-1}}{2}\,,
       \quad
   F(s) \coloneqq f(s) - f^{j-1/2}\,,
       \ \ s\in(t_{j-1},t_j), \ \ j=1,\dots,M, \\
 \intertext{and} \notag
   \psi^j \coloneqq \bigl(\m{L}U-f\bigr)^j\,,
       \ \ j=0,\dots,M,
\end{gather}
the residuum takes the form
\begin{align*}
  \left(\m{K}(u-\hat{U})\right)(s)
     = \bigl(F - \hat{F}\bigr)(s) + (t_{j-1/2}-s) \, \delta_t \psi^j + \m{L} Z^j\,.
\end{align*}
Then~\eqref{green-rep} yields
\begin{gather*}
  u(x,T) - U^M(x)
    = \sum_{j=1}^{M}
      \left\{   \int_{I_j} \scal{\m{G}(T-s)}{\bigl(F-\hat{F}\bigr)(s) + \m{L}Z^j}\D s
              + \int_{I_j} \bigl(t_{j-1/2}-s\bigr)
                \scal{\m{G}(T-s)}{\delta_t\psi^j}\D s
      \right\}.
\end{gather*}
Using~\eqref{intG} and~\eqref{intG_t-comb}, we obtain
\begin{gather*}
  \abs{\int_{I_j} \scal{\m{G}(T-s)}{\bigl(F-\hat{F}\bigr)(s)}\D s}
      \le \kappa_0 \E^{-\gamma(T-t_j)}
          \int_{I_j} \norm{\bigl(F-\hat{F}\bigr)(s)}_{\infty,\Omega} \D s\,,
\end{gather*}
\begin{gather*}
  \abs{\int_{I_j} \bigl(t_{j-1/2}-s\bigr)
                \scal{\m{G}(T-s)}{\delta_t\psi^j} \D s}
      \le \frac{\Psi_{1,j}}{2} \E^{-\gamma(T-t_j)}
          \norm{\delta_t\psi^j}_{\infty,\Omega}
\end{gather*}
and
\begin{gather}\label{extra-Z-0}
  \abs{\int_{I_j} \scal{\m{G}(T-s)}{\m{L}Z^j}\D s}
      \le \kappa_0 \tau_j \E^{-\gamma(T-t_j)}
                     \norm{\m{L}Z^j}_{\infty,\Omega}\,.
\end{gather}
Furthermore,
\begin{gather*}
  \int_{I_j} \scal{\m{G}(T-s)}{\m{L}Z^j}\D s
   = \int_{I_j} \scal{\m{L}^*\m{G}(T-s)}{Z^j}\D s
   = \int_{I_j} \scal{\pt_t\m{G}(T-s)}{Z^j}\D s
\end{gather*}
gives an alternative bound to~\eqref{extra-Z-0}:
\begin{gather*}
  \abs{\int_{I_j} \scal{\m{G}(T-s)}{\m{L}Z^j}\D s}
      \le \int_{I_j} \phi_1(T-s) \D s \
                   \E^{-\gamma(T-t_j)} \norm{Z^j(s)}_{\infty,\Omega}\,.
\end{gather*}

We arrive at the following theorem.
\begin{theorem}\label{theo:extra}
  The maximum-norm error of the extrapolated Euler
  method~\eqref{euler-extra} satisfies the a posteriori error bound
  \begin{gather*}
    \norm{u(T)-U^M}_{\infty,\Omega}
      \le \eta_{\mathrm{eE}}^M
      \coloneqq \sum_{j=1}^M \E^{-\gamma (T-t_j)}
                       \left(\eta_{\hat{F}}^j + \eta_{\delta\psi}^j + \eta^j_Z\right)
  \end{gather*}
  with the $Z^j$ defined in~\eqref{def:Zj},
  \begin{gather*}
    \eta_{\hat{F}}^j \coloneqq
       \kappa_0 \int_{I_j} \norm{\bigl(F-\hat{F}\bigr)(s)}_{\infty,\Omega} \D s\,, \quad
    \eta_{\delta\psi}^j \coloneqq
         \frac{\Psi_{1,j}}{2} \norm{\delta_t\psi^j}_{\infty,\Omega}\,, \quad
    \eta_{Z}^j \coloneqq
         \min\Biggl\{\kappa_0 \tau_j \norm{\m{L} Z^j}_{\infty,\Omega},
                    \theta_j \norm{Z^j}_{\infty,\Omega}
             \Biggr\}\,.
  \end{gather*}
\end{theorem}

\begin{remark}
  The integrals composing $\eta_{\hat{F}}$ need to be approximated.
  One possibility is Simpson's rule which gives
  \begin{gather*}
    \int_{I_j} \norm{\bigl(F-\hat{F}\bigr)(s)}_{\infty,\Omega} \D s
       \approx
       \eta_{\hat{F},\mathrm{simp}}^j
       \coloneqq  \frac{\tau_j}{6} \norm{f^j - 2 f^{j-1/2} + f^{j-1}}_{\infty,\Omega}
       \approx \frac{\tau_j^3}{24} \norm{\pt_t^2 f(t_{j-1/2})}_{\infty,\Omega}\,.
  \end{gather*}
\end{remark}

\begin{table}
\centerline{%
\begin{tabular}{c|c|cc}
$M$ & err & est & eff \\\hline
%    16 & 2.251e-04 & 1.667e-01 & 1/741 & 8.671e-02 & 4.581e-03 & 7.546e-02 \\
%    32 & 5.402e-05 & 3.506e-02 & 1/649 & 2.214e-02 & 1.356e-03 & 1.156e-02 \\
%    64 & 1.327e-05 & 7.736e-03 & 1/583 & 5.560e-03 & 3.864e-04 & 1.790e-03 \\
%   128 & 3.428e-06 & 1.798e-03 & 1/524 & 1.392e-03 & 1.052e-04 & 3.002e-04 \\
   256 & 8.780e-07 & 4.302e-04 & 1/490 \\
   512 & 2.214e-07 & 1.058e-04 & 1/478 \\
  1024 & 5.536e-08 & 2.644e-05 & 1/478 \\
  2048 & 1.382e-08 & 6.742e-06 & 1/488 \\
  4096 & 3.448e-09 & 1.699e-06 & 1/493 \\
  8192 & 8.611e-10 & 4.281e-07 & 1/497 \\
 16384 & 2.151e-10 & 1.079e-07 & 1/502 \\
 32768 & 5.369e-11 & 2.721e-08 & 1/507 \\
 65536 & 1.330e-11 & 6.859e-09 & 1/516 \\\hline
\end{tabular}}
\caption{\label{tab:extra}
  Error estimators of Theorem~\ref{theo:extra}
  for the extrapolated Euler method
  applied to test problem~\eqref{testproblem};
  Simpson's rule used to estimate $\eta_{\hat{F}}$.}
\end{table}

\begin{remark}\label{rem:asympt}
  Theorem~\ref{theo:extra} can be used to establish an asymptotically exact error
  estimator for the underlying backward-Euler discretisation:
  \begin{gather*}
    u-V = u-V + u-U = 2\left(W-V\right) + u-U.
  \end{gather*}
  Application of the triangle inequality gives
  \begin{gather*}
    \norm{u(T) - V^M}_{\infty,\Omega} \le 2\norm{W^M - V^M}_{\infty,\Omega} + \eta_{\mathrm{eE}}^M\,.
  \end{gather*}
  Similary,
  \begin{gather*}
    \norm{u(T) - W^M}_{\infty,\Omega} \le \norm{W^M - V^M}_{\infty,\Omega} + \eta_{\mathrm{eE}}^M\,.
  \end{gather*}
\end{remark}

\begin{table}
\centerline{%
\begin{tabular}{c|c|cc|cc}
M & $\norm{u(T)-W^M}_{\infty,\Omega}$ & est & eff
  & $\norm{W^M-V^M}_{\infty,\Omega}$ & $\eta_{\mathrm{eE}}^M$ \\
\hline
%    16 & 9.655e-04 & 1.679e-01 & 1/174 & 1.144e-03 & 1.667e-01 \\
%    32 & 4.504e-04 & 3.554e-02 & 1/79 & 4.892e-04 & 3.506e-02 \\
%    64 & 2.162e-04 & 7.961e-03 & 1/37 & 2.247e-04 & 7.736e-03 \\
%   128 & 1.055e-04 & 1.905e-03 & 1/18 & 1.074e-04 & 1.798e-03 \\
   256 & 5.203e-05 & 4.827e-04 & 1/9 & 5.249e-05 & 4.302e-04 \\
   512 & 2.582e-05 & 1.317e-04 & 1/5 & 2.593e-05 & 1.058e-04 \\
  1024 & 1.286e-05 & 3.933e-05 & 1/3 & 1.289e-05 & 2.644e-05 \\
  2048 & 6.417e-06 & 1.317e-05 & 1/2 & 6.424e-06 & 6.742e-06 \\
  4096 & 3.205e-06 & 4.905e-06 & 1/2 & 3.207e-06 & 1.699e-06 \\
  8192 & 1.602e-06 & 2.030e-06 & 1/1 & 1.602e-06 & 4.281e-07 \\
 16384 & 8.006e-07 & 9.087e-07 & 1/1 & 8.007e-07 & 1.079e-07 \\
 32768 & 4.003e-07 & 4.275e-07 & 1/1 & 4.003e-07 & 2.721e-08 \\
 65536 & 2.001e-07 & 2.070e-07 & 1/1 & 2.001e-07 & 6.859e-09 \\
\hline
\end{tabular}}
\caption{\label{tab:bE-asympt}
  Asymptotically exact error estimation for the backward Euler method
  according to Remark~\ref{rem:asympt}.}
\end{table}

\paragraph{Numerical results}
for the extrapolated Euler method are given in Table~\ref{tab:extra}.
They are clear illustrations for the bounds given in Theorem~\ref{theo:extra}.
The efficiency is around $500$, but slowly decreasing (with $\ln M$) as
the mesh is refined.

Table~\ref{tab:bE-asympt} illustrates Remark~\ref{rem:asympt}. Using extrapolation,
an asymptotically exact error estimator for the underlying Euler method is
obtained.
This kind of error control for initial-value problems is well established,
see, e.g., \cite[II.4]{MR1227985}: A higher-order method is used to estimate
the error of a lower-order method.
However, this approach does not guarantee upper bounds for the discretisation
error, because the error of the higher-order method is not controlled.
Additional bounds like Theorem~\ref{theo:extra} cure this defect.

\section{Discontinuous Galerkin method, dG(1)}
\label{sect:dG1}

Given $U^0=u^0$, we seek approximations $U^{j-2/3}, U^j \in H_0^1(\Omega)$ of
$u(t_{j-2/3})$ and $u(t_j)$ as solutions of
\begin{subequations}\label{dG1-discr}
\begin{align}
  U^{j-2/3} - U^{j-1} + \frac{\tau_j}{12} \left(5 \m{L} U^{j-2/3} - \m{L} U^j\right)
     & = \frac{\tau_j}{12} \left(5 f^{j-2/3} - f^j\right) \\
  U^j - U^{j-1} + \frac{\tau_j}{4} \left(3 \m{L} U^{j-2/3} + \m{L} U^j\right)
     & = \frac{\tau_j}{4} \left(3 f^{j-2/3} + f^j\right), \quad j=1,\dots,M.
\end{align}
\end{subequations}
Let $\psi\coloneqq f - \m{L}U$.
Then~\eqref{dG1-discr} can be rewritten as
\begin{gather}\label{dG1-discr*}
  U^{j-2/3} - U^{j-1} = \frac{\tau_j}{12} \left(5 \psi^{j-2/3} - \psi^j\right), \quad
  U^j - U^{j-1} = \frac{\tau_j}{4} \left(3 \psi^{j-2/3} + \psi^j\right), \quad j=1,\dots,M.
\end{gather}

We summarise the analysis from~\cite[\S6]{MR3056758} and set
\begin{gather*}
  \zeta(s) \coloneqq 3(s-1)(s-1/3) \quad \text{and} \quad
  Z(s) \coloneqq \int_0^s \zeta\left(\sigma\right) \D\sigma
               = s(s-1)^2\,,
\end{gather*}
and note that $\zeta'(s) = 6(s-2/3)$.

Given a function $v$,
we define a piecewise linear (possibly discontinuous) interpolant $\bar{v}$ by
\begin{align*}
  \bar{v}(t) & \coloneqq v^j - \frac{3}{2}\frac{t_j-t}{\tau_j} \left(v^j - v^{j-2/3}\right),
      \ \ t\in (t_{j-1},t_j], \\
  \intertext{and a continuous piecewise quadratic interpolant $\breve{v}$ by}
  \breve{v}(t) & \coloneqq v^j - \frac{3}{2}\frac{t_j-t}{\tau_j} \left(v^j - v^{j-2/3}\right)
                  + \frac{v^j - 3 v^{j-2/3} + 2 v^{j-1}}{2}
                    \zeta \left(\frac{t-t_{j-1}}{\tau_j}\right), t\in \bar{I}_j,
\end{align*}
Then, by~\eqref{dG1-discr*}
\begin{gather*}
  \breve{U}'(t)=
     \frac{3\bigl(U^j - U^{j-2/3}\bigr)}{2\tau_j}
                  + 3 \frac{U^j - 3 U^{j-2/3} + 2 U^{j-1}}{\tau_j}
                      \frac{t-t_{j-2/3}}{\tau_j}
   % = \cdots
   = \bar{\psi}(t)
\end{gather*}
This yields for the residuum
\begin{align*}
  \m{K}\bigl(u-\breve{U}\bigr)(t)
    & = f(t) - \bigl(\breve{U}' + \m{L}\breve{U}\bigr)(t)
      = \bigl(f-\breve{f}\bigr)(t) - \breve{U}'(t) + \breve{\psi}(t)
      = \bigl(f-\breve{f}\bigr)(t) + \bigl(\breve{\psi} - \bar{\psi}\bigr)(t) \\
    & = \bigl(f -\breve{f}\bigr)(t)
          + \frac{\psi^j - 3 \psi^{j-2/3} + 2 \psi^{j-1}}{2}
            \zeta \left(\frac{t-t_{j-1}}{\tau_j}\right), \ \ t\in \bar{I}_j.
\end{align*}
Set
\begin{gather}\label{def-chi-dG1}
  \chi^j \coloneqq \frac{\psi^j - 3 \psi^{j-2/3} + 2 \psi^{j-1}}{2 \tau_j^2}
     \,, \ \ j=1,\dots,M.
\end{gather}
Then the residuum can be rewritten into
\begin{align*}
  \left(\m{K}\bigl(u-\breve{U}\bigr)\right)(t)
    & = \bigl(f -\breve{f}\bigr)(t)
          + 3 \chi^j \bigl(t-t_j\bigr)\bigl(t-t_{j-2/3}\bigr) \\
    & = \bigl(f -\breve{f}\bigr)(t)
          + \chi^j \frac{\D}{\D t} \left[\bigl(t-t_j\bigr)^2
                     \bigl(t-t_{j-1}\bigr)\right], \ \ t\in \bar{I}_j,
\end{align*}
where we have used integration by parts.
Next, we multiply by the Green's function and integrate over $(0,T)$
to obtain the following a posteriori error bound.

\begin{theorem}\label{theo:dG1}
  The error of the 3rd order discontinuous Galerkin method~\eqref{dG1-discr}
  satisfies
  \begin{gather*}
    \norm{u(T)-U^M}_{\infty,\Omega}
      \sum_{j=1}^M \E^{-\gamma (T-t_j)}
                   \left( \eta_{\breve{f}}^j + \eta_{\chi}^j \right)
  \end{gather*}
  with $\chi^j$ from \eqref{def-chi-dG1},
  \begin{gather*}
    \eta_{\breve{f}}^j \coloneqq
             \int_{I_j} \norm{\bigl(f-\breve{f}\bigr)(s)}_{\infty,\Omega} \D s
      \quad\text{and}\quad
    \eta_\chi^j \coloneqq \Psi_{2,j} \norm{\chi^j}_{\infty,\Omega}\,.
  \end{gather*}
\end{theorem}

This result is a slight improvement over Theorem~6.1 in~\cite{MR3056758} as it
employs local bounds for the Green's function rather then a global argument.
An a posteriori error bound for the dG(1)-method is also given
in~\cite[\S1, Theorem 1.3]{MR1335652}, but without a proof and without
fixing the constants.
Furthermore, a remark in~\cite{MR1335652} suggests this bound is only $2$nd order
time accurate, while Theorem~\ref{theo:dG1} provides a bound of order $3$.

\begin{remark}
  Again, the integral defining $\eta_{\breve{f}}^j$ needs to be approximated.
  Simpson's rule can be applied to give
  \begin{align*}
    \int_{I_j} \norm{\bigl(f-\breve{f}\bigr)(s)}_{\infty,\Omega} \D s
       \approx \frac{2\tau_j}{3} 
        \norm{\bigl(\breve{f}-f\bigr)(t_{j-1/2})}_{\infty,\Omega}
       = \frac{2\tau_j}{3} 
        \norm{\frac{f^j+9f^{j-2/3}-2 f^{j-1}}{8} - f^{j-1/2}}_{\infty,\Omega}
    \eqqcolon f_{\breve{f},\mathrm{simp}}^j.
  \end{align*}
\end{remark}

\paragraph{Numerical results} for the dG(1)-method are presented in Table~\ref{tab:dG1}.
The results are in agreement with Theorem~\ref{theo:dG1}.
Again, looking at $M=2^{10},\dots,2^{14}$, we witness a slight deterioration
(with $\ln M$) when the mesh is refined.
For larger $M$ we are operating close to machine accuracy and the results get
erratic.

\begin{table}
\centerline{%
\begin{tabular}{c|cc|cc}
   $M$ & err & ord & est & eff \\\hline
%    16 & 1.1561e-04 & 2.60 & 2.2429e-02 & 1/194 \\
%    32 & 1.9012e-05 & 2.82 & 2.8829e-03 & 1/152 \\
%    64 & 2.6969e-06 & 2.67 & 3.6548e-04 & 1/136 \\
%   128 & 4.2436e-07 & 2.64 & 4.5565e-05 & 1/107 \\
   256 & 6.799e-08 & 2.79 & 5.739e-06 & 1/84 \\
   512 & 9.859e-09 & 2.93 & 7.270e-07 & 1/74 \\
  1024 & 1.296e-09 & 2.99 & 9.225e-08 & 1/71 \\
  2048 & 1.631e-10 & 3.00 & 1.170e-08 & 1/72 \\
  4096 & 2.032e-11 & 3.01 & 1.481e-09 & 1/73 \\
  8192 & 2.531e-12 & 3.00 & 1.872e-10 & 1/74 \\
 16384 & 3.169e-13 & 2.92 & 2.364e-11 & 1/75 \\
 32768 & 4.178e-14 & 0.66 & 2.984e-12 & 1/71 \\
 65536 & 2.645e-14 & 0.00 & 3.771e-13 & 1/14 \\
\hline
\end{tabular}}
\caption{\label{tab:dG1}
  Error estimator of Theorem~\ref{theo:dG1} for dG(1)
  applied to the test problem~\eqref{testproblem}.
}
\end{table}

\section{BDF-2}
\label{sect:BDF2}

The backward differentiation formulae (BDF-$k$) are a family multistep methods
for the approximation of initial-(boundary) value problems, and commonly used
for stiff problems.
Here we restrict ourselves to the simplest BDF-2 version, higher-order BDF-methods
are studied in~\cite{MOssad22} too.

Given $U^0=u^0$, we seek approximations $U^j \in H_0^1(\Omega)$ of
$u(t_j)$ as solutions of
\begin{subequations}\label{bdf-discr}
  \begin{align}
    \label{bdf-discr:1}
    \delta_t U^1 + \m{L} U^1 & = f^1 \\
    \label{bdf-discr:j}
   D_t U^j + \m{L} U^j & = f^j, \ \ j=2,3,\dots,M,
  \end{align}
\end{subequations}
where
\begin{gather*}
  D_t v^n \coloneqq \delta_t v^n + \tau_n \delta_t^2 v^n, \ \ \
  \delta_t^2 v^n \coloneqq \frac{\delta_t v^n - \delta_t v^{n-1}}{\tau_n + \tau_{n-1}}
  \delta_t v^n \coloneqq \frac{v^n - v^{n-1}}{\tau_n}
\end{gather*}
Again, we extend the $U^j$ to a piecewise linear function $\hat{U}$
defined on $[0,T]$.

On the first interval, the discretisation \eqref{bdf-discr:1} consists of a
single step of the implicit Euler method~\eqref{Euler_method}.
In view of our discussions following~\ref{theo:combined}, we use the argument that
led to Theorem~\ref{theo:DLMKL}.

For $s\in(t_{j-1},t_j)$, $j=2,3,\dots,M$, the residuum satisfies
\begin{align*}
  \m{K}\bigl(u-\hat{U}\bigr)(s)
    & = f(s) - \pt_t \hat{U}(s) - \m{L}\hat{U}(s) \\
    & = \bigl(f-\hat{f}\bigr)(s) -\delta_t U^j + \bigl(f-\m{L}U\bigr)^j
        + \frac{\bigl(f-\m{L}U\bigr)^j - \bigl(f-\m{L}U\bigr)^{j-1}}{\tau_j}
                      \left(s-t_j\right).
\end{align*}
By~\eqref{bdf-discr:j} we have
\begin{gather*}
  \bigl(f-\m{L}U\bigr)^j
    = \begin{cases}
        \delta_t U^j, & j=1, \\
        \delta_t U^j + \tau_j \delta_t^2 U^j, & j=2,\dots,M.
      \end{cases}
\end{gather*}
Thus,
\begin{align*}
  \m{K}\bigl(u-\hat{U}\bigr)(s)
    & = \bigl(f-\hat{f}\bigr)(s)
           + 2\left(s-t_{j-1/2}\right) \delta_t^2 U^j
           + \left(s-t_j\right) \frac{\tau_{j-1}}{\tau_j}
                \left(\delta_t^2 U^j - \delta_t^2 U^{j-1}\right)\,, \\
    & \hspace{18em} \ \ s\in(t_{j-1},t_j), \ \  j=3,\dots,M, \\
  \intertext{and}
  \m{K}\bigl(u-\hat{U}\bigr)(s)
    & = \bigl(f-\hat{f}\bigr)(s)
           + 2\left(s-t_{j-1/2}\right) \delta_t^2 U^j
           + \left(s-t_j\right) \frac{\tau_{j-1}}{\tau_j}
                \delta_t^2 U^j\,, \ \ s\in(t_1,t_2)
\end{align*}
Multiplying with $\m{G}(T-s)$ and using both~\eqref{intG} and~\eqref{intG_t-comb},
we obtain the following result.

\begin{theorem}\label{theo:BDF2}
  The maximum-norm error of the BDF-2 discretisation~\eqref{bdf-discr}
  satisfies the a posteriori bound
  \begin{align*}
    \norm{u(T)-U^M}_{\infty,\Omega}
      & \le \E^{-\gamma (T-t_1)}
                   \left( \eta_{\bar{f}}^1 + \min\left\{\eta_{\delta U}^1, \eta_{\delta \m{L} U}^1 \right\}
                   \right) \\
      & \qquad + \E^{-\gamma (T-t_2)}
                   \left( \eta_{\hat{f}}^2 + \left(\Psi_{1,2}+\kappa_0 \frac{\tau_1\tau_2}{2}\right)
                     \norm{\delta_t^2 U^2}_{\infty,\Omega}
                   \right) \\
      & \qquad + \sum_{j=3}^M \E^{-\gamma (T-t_j)}
                   \left( \eta_{\hat{f}}^j + \Psi_{1,j} \norm{\delta_t^2 U^j}_{\infty,\Omega}
                           + \kappa_0 \frac{\tau_{j-1}\tau_j}{2}
                                  \norm{\delta_t^2 U^j -\delta_t^2 U^{j-1}}_{\infty,\Omega}
                   \right)
  \end{align*}
  with $\eta_{\bar{f}}^j$ and $\eta_{\delta\m{L}U}^j$ from Theorem~\ref{theo:DLM}
  and $\eta_{\delta U}^j$ from Theorem~\ref{theo:KL1}.
\end{theorem}

\begin{remark}\label{rem:bdf-k}
  The term $\delta_t^2 U^j -\delta_t^2 U^{j-1}$ is a difference quotient
  of order $3$. For a BDF-$k$ method the technique developped in~\cite{MOssad22}
  involves difference quotients of order $2k-1$.
  Also note, that in the above analysis we had to consider the first $2$ time
  steps separately.
  For the BDF-$k$ method different arguments will be required for the first
  $2(k-1)$ steps.
\end{remark}

\begin{table}
\centerline{%
\begin{tabular}{c|c|c|c|cccc}
$M$ & err & eta & eff & \\\hline
%    16 & 3.1593e-04 & 2.3137e-01 & 1/732 \\
%    32 & 8.1301e-05 & 6.4734e-02 & 1/796 \\
%    64 & 2.0194e-05 & 1.7885e-02 & 1/886 \\
%   128 & 4.8968e-06 & 4.9795e-03 & 1/1017 \\
   256 & 1.1943e-06 & 1.4205e-03 & 1/1189 \\
   512 & 2.9418e-07 & 4.0910e-04 & 1/1391 \\
  1024 & 7.2959e-08 & 1.0916e-04 & 1/1496 \\
  2048 & 1.8164e-08 & 1.7075e-05 & 1/940 \\
  4096 & 4.5315e-09 & 4.3742e-06 & 1/965 \\
  8192 & 1.1317e-09 & 1.1241e-06 & 1/993 \\
 16384 & 2.8277e-10 & 2.9030e-07 & 1/1027 \\
 32768 & 7.0673e-11 & 7.4963e-08 & 1/1061 \\
 65536 & 1.7642e-11 & 1.9242e-08 & 1/1091 \\
\hline
\end{tabular}}
\caption{\label{tab:BDF2}
  Error estimator of Theorem~\ref{theo:BDF2} for BDF-2
  applied to the test problem~\eqref{testproblem}.
}
\end{table}

\paragraph{Numerical results} for the BDF-2 method are given in Table~\ref{tab:BDF2}.
There is a jump in the efficiency when going from $M=2^{10}$ to $M=2^{11}$
we do not have an explanation for.
Apart from this, a slight deterioration (with $\ln M$) is observed again,
when the mesh is refined.

\section{Summary and open questions}

In this paper we have reexamined (and improved) a posteriori error bounds
for semidiscretisations of parabolic PDEs. In particular we have considered
\begin{itemize*}
  \item the backward Euler method,
  \item the Crank-Nicolson method,
  \item the extrapolated Euler method
  \item the discontinuous Galerkin method with polynomial degree $1$, dG(1),
    and
  \item the BDF-2 method.
\end{itemize*}
Numerical experiments have be conducted for those methods.
They showed that the error are overestimated by a factor ranging from $50$
to $1000$. A natural question that arises is: \emph{Can these estimates be improved
to give sharper error bounds.} Ideally, one likes the efficiency of the
estimators to be close to $1$.
But there are further questions that need attention.

Richardson extrapolation: \emph{Is there an elegant way to derive error bounds
for extrapolation of arbitrary order in a common framework?}

Discontinuous Galerkin: The technique derived in~\cite[\S6]{MR3056758} for the
dG($r$) method gives a posteriori bounds with time accuracy of order $r+2$, while
the method converges with order $2r+1$.
Thus for $r\ge2$ there is a discrepancy, and the efficiency of the estimator
decays with the number of time steps (to the power of $r-1$).
\emph{Is there an alternative analysis that gives efficient a posteriori 
estimators for the dG($r$) methods?}

The backward differentiation formulae (BDF-$k$):
As noted in Remark~\ref{rem:bdf-k} the estimators derived in~\cite{MOssad22}
involve difference quotients of order $2k-1$, while $k+1$ seems to be the
natural order.
Further complications arise from the necessity to have $k$ starting values.
Again: \emph{Is there an elegant way to derive error bounds
for BDF methods of arbitrary order in a common framework?}

Continuous Galerkin: Except for the special case of Crank-Nicolson no results
are available yet.

Finally, estimators for operator splitting methods and ADI methods seem to
be desirable.

\def\cprime{$'$}

\end{document}